\documentclass[11pt]{amsart}
\usepackage{amsmath,amssymb,amsthm,mathdots,mathrsfs,mathabx}
\usepackage{graphicx} 

\usepackage[english]{babel}
\usepackage{algorithm}
\usepackage[noend]{algpseudocode}

\theoremstyle{plain}
\newtheorem{theorem}{Theorem}[section]
\newtheorem{lemma}[theorem]{Lemma}

\theoremstyle{definition}
\newtheorem{remark}[theorem]{Remark}
\newtheorem{definition}[theorem]{Definition}
\newtheorem{example}[theorem]{Example}

\newtheorem{problem}[theorem]{Problem}

\numberwithin{equation}{section}

\newcommand{\BC}{{\mathbb C}}

\newcommand{\BH}{{\mathbb H}}

\newcommand{\BN}{{\mathbb N}}
\newcommand{\BP}{{\mathbb P}}
\newcommand{\BR}{{\mathbb R}}

\newcommand{\cL}{{\mathcal L}}
\newcommand{\cM}{{\mathcal M}}
\newcommand{\cO}{{\mathcal O}}\newcommand{\cP}{{\mathcal P}}
\newcommand{\cR}{{\mathcal R}}

\newcommand{\bS}{{\mathbf S}}

\newcommand{\wtilA}{\widetilde{A}}\newcommand{\wtilB}{\widetilde{B}}

\newcommand{\wtilE}{\widetilde{E}}

\newcommand{\Ga}{\Gamma}

\newcommand{\la}{\lambda}

\newcommand{\om}{\omega}
\newcommand{\rank}{\operatorname{rank}}
\newcommand{\ran}{\operatorname{Ran}}

\newcommand{\im}{\operatorname{Im}}
\newcommand{\re}{\operatorname{Re}}
\newcommand{\kr}{\operatorname{Ker}}
\newcommand{\diag}{\operatorname{diag}}

\newcommand{\ov}[1]{{\overline{#1}}}

\newcommand{\pmat}[1]{\begin{pmatrix} #1 \end{pmatrix}}

\newcommand{\sbm}[1]{\left(\begin{smallmatrix} #1\end{smallmatrix}\right)}

\newcommand{\PRO}{\mathcal{PRO}}
\newcommand{\PO}{\mathcal{PO}}

\newcommand{\half}{\frac{1}{2}}

\newcommand{\Ree}{\operatorname{Re}}
\newcommand{\Sym}{\operatorname{Sym}}

\title[Interpolation for positive real odd functions]{Matrix-point interpolation for positive real odd functions in one and several variables}

\author[S. ter Horst]{Sanne ter Horst}
\address{S. ter Horst, School of Mathematical and Statistical Sciences,
North-West University,
Research Focus:\ Pure and Applied Analytics,
Private Bag X6001,
Potchefstroom 2520,
South Africa
and
DSI-NRF Centre of Excellence in Mathematical and Statistical Sciences (CoE-MaSS),
Johannesburg,
South Africa}
\email{Sanne.TerHorst@nwu.ac.za}

\author[A. van der Merwe]{Alma van der Merwe}
\address{A. van der Merwe, Department of Mathematics and Applied Mathematics, University of Pretoria, Private bag X20 Hatfield, 0028 Pretoria, South Africa}
\email{alma.vandermerwe@up.ac.za}

\author[V. Vinnikov]{Victor Vinnikov}
\address{V. Vinnikov, Department of Mathematics, Ben-Gurion University of the Negev, P.O. Box 653, Beer-Sheva 8410501, Israel}
\email{vinnikov@math.bgu.ac.il}

\author[H.J. Woerdeman]{Hugo J. Woerdeman}
\address{H.J. Woerdeman, Department of Mathematics, Drexel University, Philadelphia, PA 19104, USA}
\email{hugo@math.drexel.edu}

\thanks{This work is based on the research supported in part by the National Research Foundation of South Africa (Grant Number 127364). In addition, HW is partially supported by NSF grant DMS-2348720.}



\begin{document}

\maketitle

\begin{abstract}
We study the single matrix point interpolation problem for positive real matrix valued functions on the right half-plane, both in the single variable setting and in the noncommutative multivariable setting. In the single variable case, a careful analysis is done for when the solution may or must have a pole at infinity. In the multivariable setting, as the pole structure at infinity can be more complicated, we focus on the case when there is no pole at infinity. The results also yield numerical algorithms for finding a solution.
\end{abstract}

\section{Introduction}

We write $\BC(z)$ and $\BR(z)$ for the rational functions in the complex variable $z$ with complex and real coefficients, respectively. We write $\cP^{m\times m}$ for the class of {\em positive real} $m\times m$ rational matrix functions $f$ without poles in the open right half-plane ${\mathbb C}_+:=\{z\in\BC \colon \re z >0\}$, i.e., such that \[
\re f(z)\geq 0 \mbox{ whenever } \re z>0. \]
Next, $\cO^{m\times m}$ denotes the class of {\em odd} $m\times m$ rational matrix functions $f$, i.e.,
\[ f(z)=-{f(-\ov{z})}^*. \]
In addition, we define the class of {\em positive real odd matrix functions} $\PO^{m\times m}:=\cP^{m\times m}\cap \cO^{m\times m}$.
If the superscript $m\times m$ is missing, then $m=1$.
We also have real variants of these classes of functions, such as $\cP\cR:=\cP\cap \BR(z)$, $\PRO:=\PO\cap \BR(z)$, etc. Scalar valued positive real odd functions were first introduced in the electrical engineering literature (see Foster \cite{F24}, Cauer \cite{C26} and Brune \cite{B31}; see also \cite{N66,B68,AV73}).

A function $f\in\BC(z)$, written as $f(z)=\frac{p(z)}{q(z)}$ for co-prime polynomials $p$ and $q$ is evaluated in a matrix $A\in \BC^{n \times n}$ via $f(A)=p(A)q(A)^{-1}$ whenever $q(A)^{-1}$ exists, that is, whenever the intersection of the eigenvalues of $A$ and the roots of $q$ is empty. By standard functional calculus, the definitions of the classes $\cP^{m\times m}$ and $\cO^{m\times m}$ extend to matrix-point evaluation. Hence $f\in \cP^{m\times m}$ if and only if for each $n\in\BN$ and each $A\in\BC^{n \times n}$ with $\re A:=\half(A+A^*)>0$, we have $\re f(A)\geq 0$, whereas $f\in \cO^{m\times m}$ holds if and only if for each $n\in\BN$ and $A\in\BC^{n \times n}$  we have $-f(A)=f(-A^*)^*$. 

The main problem in this paper concerns the single matrix point interpolation problem, which in the case of one variable is the following. 

\begin{problem}\label{Problem}
For $A\in \BC^{n \times n}$ and $B\in \BC^{nm \times nm}$, when does there exist a function $f$ in $\cP^{m\times m}$, $\PO^{m\times m}$ or $\PRO^{m\times m}$ such that $f(A)=B$? If so, provide a way to find a/all solutions.
\end{problem}

Basic functional calculus shows that $B$ should be in the bicommutant of $A$, and in that case, when $A$ and $B$ are also Lyapunov regular, one can reduce to the case where $A$ is nonderogatory and $\re A$ positive definite. Indeed, in the Jordan canonical form of $A$ one can simply remove  any blocks with repeated eigenvalues, leaving only the largest one, and reflect eigenvalues of $A$ with negative real part over the imaginary axis, making appropriate adjustments in $B$.

This problem was first considered in the real scalar case (thus, $\PRO$) by N. Cohen and I. Lewkowicz (in \cite{CL07,CL09}), for the case that $A$ and $B$ are Lyapunov regular, that is, in both cases the eigenvalues $\la_j$ satisfy $\la_j+\ov{\la}_i\neq 0$ for all $i$ and $j$. More recently, S. ter Horst and A. van der Merwe
\cite{tHvdM22b} studied the problem in $\PO$ and $\PRO$ and identified the so-called Hill-Pick matrix $\BH_{A,B}$ of the pair $(A,B)$ and showed that solvability of the interpolation problem requires $\BH_{A,B}$ to be positive semidefinite. With the additional condition that  $\rank \BH_{A,B}$ equals the dimension of the real bicommutant of $A$, they showed in \cite{tHvdM} that a solution exists and provided a constructive way to find the solution.  This solution $f(z)$ has the additional property $\lim_{z\to\infty} f(z)=0$, which is referred to as the solution having no pole at infinity. 

In what follows we build on the results by Ter Horst and Van der Merwe, and generalize them in several directions. Namely, in the matrix case, we prove existence of a solution without the condition on $\rank \BH_{A,B}$, in which a pole at $\infty$ can occur, and we clarify in detail when that happens, as well as when there exists a unique solution. We then extend some of the results to the multivariable noncommutative function setting. Here, with a domain described via a linear matrix inequality $Q(z)>0$, we first establish necessary and sufficient condition for the existence of an interpolant with a realization of the form \begin{equation}\label{fzDM} f(z)=
D
+\ell^*
\left(I_s\otimes Q(z)-M\right)^{-1}
\ell,\ D=-D^*, M=-M^*. \end{equation} These conditions involve complete positivity of a map defined on an operator system and the existence of an Arveson extension, where the extension is required to satisfy an additional constraint; see Theorem \ref{thm:Q-interpolation} for details. 
Next, we consider the special case where $Q(z) = {\rm diag}(z_j)_{j=1}^d$; that is, the case when the domain is the standard noncommutative poly-half-plane. In this case, we can give a functional theoretic characterization of rational matrix functions in noncommutative variables of the form \eqref{fzDM}; see Theorem \ref{POd}. We end the paper with providing a computational algorithm (using SemiDefinite Programming) to find a multivariable interpolant of the form \eqref{fzDM} and illustrate the results in some examples. 

The paper is organized as follows. 
In Section 2 we will discuss in detail the one variable case, where already many results exist. 
Putting the matrix $A$ in Jordan canonical form allows for detailed computational descriptions.
In Section 3, we provide a view of the one variable case that includes completely positive maps, which can then be used to develop partial results in the multivariable noncommutative case, which we will treat in Sections 4 (generalized poly-half planes) and 5 (standard poly-half plane). The approach here is a more coordinate free one, in part as there is no convenient canonical form for a tuple of noncommuting matrices.
Along the way, we illustrate the results with examples and provide constructions to find a solution.

\section{The one variable case}


In this section, we discuss the one variable interpolation results. The results are treated in different parts of the literature, so this section serves as a convenient recap. As most of the results are not new, we do not provide all details in the proofs but rather provide the main ideas. In the next section we will approach the one variable case from a different and coordinate free perspective, as a lead into the multivariable case. The spectrum of a square matrix $M$ is denoted by $\sigma(M)$.

\subsection{Background on state space realizations}
More details and proofs on the results we recall in this subsection can be found in \cite{Kailath1980,BGR, Dai1989,KunkelMehrmann2006}.
We let rational functions over $\mathbb F\in\{\mathbb C,\mathbb R\}$  be denoted by \(\mathbb{F}(z)\), which is the field of fractions of the polynomial ring \(\mathbb{F}[z]\) of one variable.
\begin{definition}[Minimal realization]\label{D:MinReal}
Let
\[
G(z)\in\mathbb{F}(z)^{p\times m}
\]
be a matrix-valued rational function over a field $\mathbb{F}$. A \emph{state-space realization} of $G$ is a quadruple
\[
(A,B,C,D),
\]
with
\[
A\in\mathbb{F}^{n\times n},\qquad
B\in\mathbb{F}^{n\times m},\qquad
C\in\mathbb{F}^{p\times n},\qquad
D\in\mathbb{F}^{p\times m},
\]
such that
\[
G(z)=D+C(zI-A)^{-1}B.
\]

The realization $(A,B,C,D)$ is called \emph{minimal} if its state dimension $n$ is the smallest among all realizations of $G(z)$; that is, if there does not exist another realization
\[
(\widetilde A,\widetilde B,\widetilde C,\widetilde D)
\]
with state dimension $\widetilde n<n$ satisfying
\[
G(s)=\widetilde D+\widetilde C(sI-\widetilde A)^{-1}\widetilde B.
\]
The minimal state dimension is called the \emph{McMillan degree} of $G$.
\end{definition}
\begin{theorem}
A realization $(A,B,C,D)$ of the rational matrix $G(s)$ is minimal if and only if it is both controllable and observable. Equivalently,
\[
\rank
\begin{pmatrix}
B & AB & \cdots & A^{n-1}B
\end{pmatrix}
=n,
\]
and
\[
\rank
\begin{pmatrix}
C\\
CA\\
\vdots\\
CA^{n-1}
\end{pmatrix}
=n.
\]
\end{theorem}

Rational matrix function given by a realization as in Definition \ref{D:MinReal} have the limitation that they cannot have a pole at $\infty$. Since the functions we consider in Problem \ref{Problem} can have poles at $\infty$, we will also consider  descriptor systems which we define next; the following can be found in \cite[Sections 2-2, 2-3]{Dai1989}.

\begin{definition}[Minimal descriptor realization]
Let
\[
G(z)\in\mathbb{F}(z)^{p\times m}
\]
be a matrix-valued rational function. A \emph{descriptor realization} (or \emph{singular realization}) of $G$ is a quintuple
\[
(E,A,B,C,D),
\]
where
\[
E,A\in\mathbb{F}^{n\times n},\qquad
B\in\mathbb{F}^{n\times m},\qquad
C\in\mathbb{F}^{p\times n},\qquad
D\in\mathbb{F}^{p\times m},
\]
such that the matrix pencil $zE-A$ is regular, i.e.,
\[
\det(zE-A)\not\equiv 0,
\]
and
\[
G(z)=D+C(zE-A)^{-1}B.
\]

The descriptor realization $(E,A,B,C,D)$ is called \emph{minimal} if its state dimension $n$ is the smallest among all regular descriptor realizations of $G(z)$. The minimal state dimension is called the McMillan degree of $G$.
\end{definition}

\begin{definition}
Consider the descriptor realization $(E,A,B,C,D)$ of a rational matrix function $G(z)\in\mathbb{F}(z)^{p\times m}$.

\begin{enumerate}
    \item The triple $(E,A,B)$ is said to be \emph{finite controllable} if
    \[
    \rank
    \begin{pmatrix}
    \lambda E-A & B
    \end{pmatrix}
    = n
    \]
    for every finite generalized eigenvalue $\lambda$ of the pencil $sE-A$.

    \item The triple $(E,A,B)$ is said to be \emph{impulse controllable} if
    \[
    \rank
    \begin{pmatrix}
    E & 0 & 0 \\
    A & E & B
    \end{pmatrix}
    = n +\rank \ E.
    \]

    \item The triple $(E,A,C)$ is said to be \emph{finite observable} if
    \[
    \rank
    \begin{pmatrix}
    \lambda E-A\\
    C
    \end{pmatrix}
    = n
    \]
    for every finite generalized eigenvalue $\lambda$ of the pencil $sE-A$.

    \item The triple $(E,A,C)$ is said to be \emph{impulse observable} if
        \[
    \rank
    \begin{pmatrix}
    E  & A \\ 0 & E\\
     0 & C
    \end{pmatrix}
    = n + \rank E.
    \]
\end{enumerate}
\end{definition}

\begin{theorem}\label{T:DescrRealMin}
Let $(E,A,B,C,D)$ be a descriptor realization of a rational matrix function $G(z)\in\mathbb{F}(z)^{p\times m}$.
Then the following are equivalent:
\begin{enumerate}
    \item $(E,A,B,C,D)$ is a minimal descriptor realization;
    \item
    \begin{enumerate}
        \item $(E,A,B)$ is finite controllable;
        \item $(E,A,B)$ is impulse controllable;
        \item $(E,A,C)$ is finite observable;
        \item $(E,A,C)$ is impulse observable.
    \end{enumerate}
\end{enumerate}
\end{theorem}


\subsection{Realization theory for functions in $\cP^{m\times m}$, $\cP\cR^{m\times m}$, $\cP\cO^{m\times m}$ and $\cP\cR\cO^{m\times m}$}

In this subsection we provide state space realization formulas for functions in the classes $\cP^{m\times m}$, $\cP\cO^{m\times m}$ and $\cP\cR\cO^{m\times m}$, starting with $\cP^{m \times m}$ and $\cP\cR^{m\times m}$.

\begin{theorem}\label{P}
Let $f\in\BC(z)^{m \times m}$. Then $f \in \cP^{m\times m}$ if and only if $f$ can be written as 
\begin{equation}\label{0}     
f(z) = Kz + D+ C(zI_k+A)^{-1}B,
\end{equation}
where $k\in {\mathbb N}$, $K\in \BC^{m\times m}$ and the block matrix $\pmat{A&B\\C&D} \in\BC^{(k+m)\times(k+m)}$ are such that
\begin{equation}\label{eq:c2}
K \ge 0,\quad \pmat{
A+ A^* & C^*- B \\ C- B & D+D^*
} \geq 0.
\end{equation}
The above realization may be chosen to be minimal, meaning that there exists no representation of $f$ of the form \eqref{0} with state matrix $A$ of smaller size.

Moreover, $f \in \cP\cR^{m\times m}$ if and only if the matrices $A,B,C,D$, and $K$ may be chosen to be real.
\end{theorem}

\begin{proof} 
For $\cP\cR^{m\times m}$, the result follows from Theorems 1 and 3 in \cite{Willems}. The arguments for $\cP^{m\times m}$ are analogous.
\end{proof}

Next we provide realization formulas for functions in the classes $\PO^{m\times m}$ and $\PRO^{m\times m}$.

\begin{theorem}\label{PO} For $f\in\BC(z)^{m \times m}$ the following are equivalent. \begin{itemize} 
\item[(i)] $f \in \PO^{m\times m}$; 
    
\item[(ii)] $f$ may be written as 
\begin{equation}\label{1}     f(z) = Kz + D+ C(zI_k+A)^{-1}B,
\end{equation}
where $k\in {\mathbb N}$, $K\in \BC^{m\times m}$ and the block matrix $\pmat{A&B\\C&D} \in\BC^{(k+m)\times(k+m)}$ are such that
\begin{equation}\label{eq:c2}
K \ge 0,\quad D=-D^*,\quad B=C^*,\quad A=-A^*;
\end{equation}

\item[(iii)]
$f$ may be written as 
    \begin{equation}\label{F}
   f(z)= Kz + D+ B^*(zI_k -i\Lambda)^{-1}B, \end{equation} 
    with $D=-D^*, K \ge 0$ and $\Lambda$ a diagonal real matrix;
    
\item[(iv)] $f$ may be written as 
\begin{equation}\label{2}     f(z) =  D+ C(zE+A)^{-1}B,
\end{equation}
where 
\begin{equation}\label{eq:c3}
E\ge 0,\   D=-D^*,\  B=C^*,\  A=-A^*,\  \kr E \cap \kr A = (0).
\end{equation}
\end{itemize}
The realizations in \eqref{1}, \eqref{F} and \eqref{2} may be chosen to be minimal, i.e., one cannot find a realization of $f$ of the forms \eqref{1}, \eqref{F} and \eqref{2} with matrices $A$ and $\Lambda$ of smaller size.

Moreover, $f(z) \in \PRO^{m\times m}$, if and only if the matrices $A,B,C,D, E$, and $K$ in \eqref{1} and \eqref{2} may be chosen to be real.
\end{theorem}

\begin{proof}
For $\PRO^{m\times m}$, the equivalence of (i) and (ii) is compiled from various sources in \cite[Theorem 3.1]{tHN20}; for $\PO^{m\times m}$  similar arguments apply. The equivalence of (ii) and (iii) follows by simply diagonalizing $A=-A^*$ via a unitary similarity. 

To see that (ii) implies (iv), let $f$ be given as in \eqref{1}. Take a minimal rank factorization $K=L^*L$, with $r=\rank K$, and note that 
\begin{align*}
f(z) &=  L^*Lz + D+ C^*(zI_r -i A)^{-1}C\\
&=D+ \pmat{0\\L\\C}^* \left(z\pmat{I_r&0&0\\0&0&0\\0&0&I_k} - \pmat{0&I_r&0\\ -I_r & 0 & 0\\ 0 & 0 & A}  \right)^{-1}
\pmat{0\\L\\C},
\end{align*}
which provides the descriptor realization of the desired form. Assume in addition that \eqref{1} is a minimal realization, so that, by the Hautus test, $\rank \pmat{z I_k -A & B}=k$ and $\rank \pmat{z I_k -A \cr C}=k$ for all $z\in\BC$. One then easily checks that the above regular descriptor realization of $f$ satisfies all conditions in item (2) of Theorem \ref{T:DescrRealMin}, to conclude that the regular descriptor realization of $f$ is minimal.   

For the converse implication, assume $f$ is given by a descriptor realization as in (iv). Recall that the matrix pencil $zE+A$ is $*$-congruent to a matrix pencil $z \wtilE +\wtilA$ when there exists an invertible matrix $S$ such that $\wtilE=S^*ES$ and $\wtilA=S^*AS$. Applying the canonical form result of \cite[Theorem 4.1(b)]{LR} to the selfadjoint pencil $E+\lambda (-iA)$, noting that certain blocks cannot appear due to $E$ being positive semidefinite, it follows that $z \wtilE +\wtilA$ is $*$-congruent to a matrix pencil that is formed by direct sum of blocks of the following types:
\begin{enumerate}
\item[\textnormal{(i)}] Finite eigenvalue blocks:
\[
\left( [1],\; [i\lambda] \right), \quad \lambda \in \mathbb{R},\quad \mbox{corresponding to scalar pencils}\quad z - i\lambda.
\]

\item[\textnormal{(ii)}] Infinite eigenvalue blocks:
\[
\left( [0],\; [i\mu] \right), \qquad \mu \in \mathbb{R},\quad \mbox{corresponding to scalar pencils}.
\]

\item[\textnormal{(iii)}] Singular (Kronecker) blocks: 
\[
\left(\begin{pmatrix}
1 & 0 \\
0 & 0
\end{pmatrix}, \begin{pmatrix}
0 & 1 \\
-1 & 0
\end{pmatrix}\right),\quad \mbox{corresponding to pencils}\quad 
\begin{pmatrix}
z & 1 \\
-1 & 0
\end{pmatrix}.
\]
\end{enumerate}
Hence $z E + A$ is $*$-congruent to a pencil $z \wtilE +\wtilA$ of the form
\begin{equation}\label{EA}
\wtilE=\pmat{I_{r_1}&0&0&0\\0&0&0&0\\0&0&I_{r_2}&0\\0&0&0&0},\quad 
\wtilA=\pmat{0&I_{r_1}&0&0\\-I_{r_1}&0&0&0\\0&0&i \Lambda_2&0\\0&0&0&i \Lambda_3}
\end{equation}
with $\Lambda_2$ and $\Lambda_3$ real diagonal matrices, where $\Lambda_3$ is invertible, since the pencil is regular. If $S$ is the invertible matrix establishing the $*$-congruence, then $(z E-A)^{-1}=S(z \wtilE - \wtilA)^{-1}S^*$, so that $f(z)=D+\wtilB^*(z \wtilE - \wtilA)^{-1}\wtilB$ with $\wtilB=S^*B$. Decompose $\wtilB$ according to the decompositions of $\wtilE$ and $\wtilA$ as
\[
\wtilB^*=\pmat{\wtilB_1^*& \wtilB_1'^* & \wtilB_2^* & \wtilB_3^*}.
\]
Then 
\begin{align*}
f(z)&=D+\wtilB^*(z \wtilE - \wtilA)^{-1}\wtilB\\ 
&=D+ i \wtilB_3^* \Lambda_3^{-1}\wtilB_3 + \wtilB_1'^*\wtilB_1- \wtilB_1^*\wtilB_1' + z \wtilB_1'^* \wtilB_1'
+ \wtilB_2^*(z I_{r_3}+i\Lambda_2)^{-1}\wtilB_2.
\end{align*}
Hence $f$ has a state space realization as in (iii). 

In (ii) and (iii), if the state space realizations are not minimal, one can apply the Kalman decomposition to make it minimal, and it is easy to check that the symmetry assumptions on the matrices are preserved in the minimal realization obtained from the Kalman decomposition. We already noted in the proof of (ii) $\Rightarrow$ (iv) that minimality is preserved, so that we can also get a minimal realization as in (iv). That the matrices $A,B,C,D$ in (ii) can be taken real if $f\in \PRO^{n \times n}$ follows from \cite[Theorem 3.1]{tHN20}. That this is also the case in (iv) follows by verifying that the matrices constructed in the implication (ii) $\Rightarrow$ (iv) remain real when $A,B,C,D$ in (ii) were real matrices. 
\end{proof}

\subsection{Solution to Problem \ref{Problem}}

Next we provide a solution to Problem \ref{Problem} for the case that $\re A>0$ and $\re B>0$. By standard functional calculus, it follows that for whenever $f(A)=B$ for a $f\in\BC(z)^{m \times m}$, we must have that $B$ is in $\{I_m \otimes A\}''$, the double commutant of $I_m \otimes A$. 

Let $A=S(J_{n_1}(\lambda_1) \oplus \cdots \oplus J_{n_l}(\lambda_l ))S^{-1}$ be the Jordan canonical form of $A$, with $S\in\BC^{n \times n}$ invertible, and $J_k(\lambda)$ is the $k\times k$ upper triangular Jordan block with eigenvalue $\lambda$. Hence, $n=n_1+\cdots + n_l$ and $\re \la_j>0$ for $j=1,\ldots,l$. Introduce the vector 
\begin{equation}\label{v}
v:=S^{*-1}\pmat{e_1^{(n_1)} \cr \vdots \cr e_1^{(n_l)}}\in \BC^n,
\end{equation}
where $e_1^{(k)}$ is the first standard basis vector of ${\mathbb C}^k$.

We can now formulate the main result of this section.

\begin{theorem}\label{ot45} 
Let $A\in\BC^{n\times n},B\in\BC^{nm\times nm}$ with $\re A >0$, $\re B >0$ and $B\in\{I_m \otimes A\}''$. Define $v\in\BC^n$ by \eqref{v}. Then the following are equivalent:
\begin{itemize}
    \item[(i)] there exists $f\in \cP^{m \times m}$ with $f(A)=B$;
    \item[(ii)] there exists $f\in \PO^{m \times m}$ with $f(A)=B$;
    \item[(iii)] $H\ge 0$, where $H$ is the (unique) solution to
    \begin{equation}\label{Lya}
(I_m\otimes A)^* H + H (I_m\otimes A) = B^*(I_m\otimes vv^*)+ (I_m\otimes vv^*)B.        
    \end{equation} 
\end{itemize}
If one (and thus all) of (i)-(iii) holds, a particular solution in $ \PO^{m \times m}$ is given by
\begin{equation}\label{fz}    
f(z)= (I_m\otimes v^*)B[zH-(H(I_m\otimes A)-B^*(I_m\otimes vv^*))]^{(-1)} B^*(I_m\otimes v).
\end{equation}
Here $T^{(-1)}$ indicates the Moore-Penrose inverse of the matrix $T$.
\end{theorem}

In the proof of Theorem \ref{ot45}, we use the following lemmas. 

\begin{lemma}\label{vprop} 
With $A, B$ and $v$ as in Theorem \ref{ot45} and $f\in\BC(z)^{m \times m}$, we have that
$f(A)=B$ if and only if $(I_m\otimes v^*)f(A)=(I_m\otimes v^*)B$.
\end{lemma}

\begin{proof}
The only if part is trivial, so suppose that $(I_m\otimes v^*)f(A)=(I_m\otimes v^*)B$. Write $f=[f_{ij}]_{i,j=1}^m$ and $B=[B_{ij}]_{i,j=1}^m$, with $B_{ij}\in\BC^{n \times n}$. It suffices to show $f_{ij}(A)=B_{ij}$ for all $i,j=1,\ldots,m$. We now have that
\begin{equation}\label{vB} \begin{pmatrix} e_1^{(n_1)*} & \dots & e_1^{(n_l)*} \end{pmatrix} (S^{-1} f_{ij}(A) S) = \begin{pmatrix} e_1^{(n_1)*} & \dots & e_1^{(n_l)*} \end{pmatrix} S^{-1} B_{ij} S. \end{equation}
Since $B_{ij} \in \{A\}''$, we have that $S^{-1} B_{ij} S$ has the same block structure as 
$$ S^{-1} A S = J_{n_1}(\lambda_1) \oplus \cdots \oplus J_{n_l}(\lambda_l),$$
with the blocks in $B_{ij}$ being upper triangular Toeplitz matrices, and for repeated eigenvalues in $A$ the corresponding blocks overlap in values. But then it suffices to read of the top row of each block (due to the Toeplitz structure), which is exactly what the left multiplication with $e_1^*$ does. This now gives that $f_{ij}(A)=B_{ij}$, and thus the lemma follows.
\end{proof}

\begin{lemma}\label{deriv} 
Let $E\ge 0$ and $A=-A^*$ and consider the pencil $zE+A$. Then
\begin{equation}\label{deriv2}
\frac{d^k}{dz^k} (zE+A)^{(-1)}= (-1)^k k! ((zE+A)^{(-1)} E)^{k}(zE+A)^{(-1)}.  
\end{equation} 
\end{lemma}

\begin{proof}
   As explained in the proof of Theorem \ref{PO}, we have that $zE+A=S^{*-1}(z\wtilE + \wtilA)S^{-1}$, where $\wtilE, \wtilA$ are as in \eqref{EA}. Then 
   $$ (zE+A)^{(-1)} = S 
\pmat{0&-I_{r_1}&0&0\\I_{r_1}&zI_{r_1}&0&0\\0&0& (zI_{r_2} + i \Lambda_2)^{-1}&0\\0&0&0&i \Lambda_3^{(-1)}} 
S^*, $$
and
\begin{align*}
\frac{d}{dz} (zE+A)^{(-1)} 
&= S 
\pmat{0&0&0&0\\0&I_{r_1}&0&0\\0&0& -(zI_{r_2} + i \Lambda_2)^{-2}&0\\0&0&0&0} S^*\\
&= -(zE+A)^{(-1)} E (zE+A)^{(-1)},
\end{align*}
which proves \eqref{deriv2} for $k=1$.  
Proceeding by induction, let us take the derivative of \eqref{deriv2}, leading to
\begin{align*}
&\frac{d^{k+1}}{dz^{k+1}} (zE+A)^{(-1)}=(-1)^k k! \left( \sum_{j=1}^{k} ((zE+A)^{(-1)} E)^{j-1} \times\right.\\
&\quad  \times [-(zE+A)^{(-1)} E (zE+A)^{(-1)}]E((zE+A)^{(-1)} E)^{k-j}(zE+A)^{(-1)} +\\
&\qquad\ \ + \left.((zE+A)^{(-1)} E)^{k}[-(zE+A)^{(-1)} E (zE+A)^{(-1)}] \right)\\
&\qquad\qquad= (-1)^{k+1} k!(k+1) ((zE+A)^{(-1)} E)^{k+1}(zE+A)^{(-1)}.\qedhere
\end{align*}
\end{proof}

In the proof below and later on, we will also make use of the {\em Lyapunov equation}
$$ XA+A^*X= H.$$ When ${\rm Re}A >0 $ the equation has a unique solution $X$ for every right hand side $H$ (see, e.g., \cite[Section 13.1]{LT}); In fact, $X=\int_0^{\infty} e^{-tA^*}He^{-tA}dt $ is the unique solution. Note that $X\ge 0$ when $H\ge 0$.

\begin{proof}[Proof of Theorem \ref{ot45}]
The implication (ii) $\Rightarrow$ (i) is trivial.

(i) $\Rightarrow$ (iii): Let $f\in \cP^{m\times m}$ such that $f(A)=B$. Let $T\in {\mathbb C}^{n\times n}$ be invertible. If $\re T AT^{-1} >0$, then due to $f\in \cP^{m\times m}$, we have that 
$\re f(TAT^{-1}) = \re\, (I_m \otimes T)B (I_m \otimes T^{-1}) \ge 0.$ Letting $P=T^*T$, this translates into
$$ P>0,\ PA+A^*P >0 \implies (I_m \otimes P)B +B^* (I_m \otimes P)  \ge 0.$$ If $P \ge 0$ and $PA+A^*P \ge 0$, one can replace $P$ by $P+\epsilon I$, apply the above, and let $\epsilon$ go to 0. This argument yields that
\begin{equation}\label{imply}
P\ge 0,\ PA+A^*P \ge 0 \implies (I_m \otimes P)B +B^* (I_m \otimes P)  \ge 0.\end{equation} 
The implication \eqref{imply} can be used to introduce the notion of Lyapunov order; see \cite{CL09} for the scalar case.
Take $P = \int_0^\infty e^{-tA^*}vv^*e^{-tA} dt \ge 0$. Then $PA+A^*P = vv^* \ge 0$. 
By  \eqref{imply}, we have $(I_m \otimes P)B +B^* (I_m \otimes P)  \ge 0.$ It remains to show that $H:=(I_m \otimes P)B +B^* (I_m \otimes P)$ satisfies \eqref{Lya}. Plugging this expression into the left hand side of \eqref{Lya}, and using that $I_m\otimes A$ and $B$ commute, we get 
\begin{align*}
(I_m\otimes A)^* H + H(I_m\otimes A)&=(I_m\otimes A)^* [(I_m \otimes P)B +B^* (I_m \otimes P)] +\\  
&\qquad\qquad +[(I_m \otimes P)B +B^* (I_m \otimes P)](I_m\otimes A)\\
&= (I_m\otimes (vv^*-PA)) B + B^* (I_m\otimes A^*P) + \\
&\qquad\qquad + (I_m\otimes PA) B + B^* (I_m\otimes (vv^*-A^*P))\\ 
&= B^*(I_m\otimes vv^*)+ (I_m\otimes vv^*)B.
\end{align*}
Thus, by the uniqueness of the solution to \eqref{Lya}, we have that $H$ satisfying \eqref{Lya} is positive semidefinite. 

(iii) $\Rightarrow$ (ii): Given (iii), to prove (ii) we show that $f(z)$ in \eqref{fz} is in $\PO^{m \times m}$ and satisfies $f(A)=B$. When $H$ is invertible, then the Moore-Penrose inverse can be replaced with a regular inverse, it follows from Theorem \ref{PO} that $f \in \PO^{m\times m}$, and by applying \cite[Theorem 22.2.2]{BGR} (with $g=0$) to the problem $(I\otimes v^*)f(A)= (I\otimes v^*)B$ together with Lemma \ref{vprop} it follows that $f(A)=B$. 

For the case where $H$ is not invertible, we provide a direct proof.  Put $\Gamma= H(I_m\otimes A)-B^*(I_m\otimes vv^*)$. It follows from \eqref{Lya} that $\Gamma = - \Gamma^*$. Let $x \in \kr H \cap \kr \Gamma$. Then $x\in \kr \,(I_m\otimes vv^*)B$ by
\begin{align*}
0
&=\Ga x = -\Ga^* x = -(I_m\otimes A^*)Hx+(I_m\otimes vv^*)Bx
=(I_m\otimes vv^*)Bx. 
\end{align*}
This implies that in \eqref{fz}, the Moore-Penrose inverse can be replaced with the regular inverse of the compression of $zH-\Gamma$ to $(\kr H \cap \kr \Gamma)^\perp$, and it then follows from Theorem \ref{PO} that $f \in \PO^{m\times m}$. 

It remains to show that $f(A)=B$. As in the proof of Lemma \ref{vprop}, write $f=[f_{ij}]_{i,j=1}^m$ and $B=[B_{ij}]_{i,j=1}^m$, so that the problem reduces to showing that $f_{ij}(A)=B_{ij}$ for $i,j=1,\ldots,m$. Also, after appropriate transformations using the similarity matrix $S$ in the Jordan canonical form of $A$ we may assume that $A$ is of the form
$$ A=J_{n_1}(\lambda_1) \oplus \cdots \oplus J_{n_l}(\lambda_l ),\quad \re\la_j>0,\ \ j=1,\ldots,l.$$
Since $B\in\{I_m\otimes A\}''$, it follows that $B_{ij}\in \{A\}''$, which implies that 
\[
B_{ij} = T_1 \oplus \cdots \oplus  T_l,
\] 
where each $T_s$ is an upper triangular Toeplitz matrix in $\BC^{n_s \times n_s}$. Hence, $f_{ij}(A)=B_{ij}$ if and only if $f_{ij}(J_{n_s}(\lambda_s ))=T_s$ for $s=1,\ldots,l$. Due to the particular form of $H$, the problem reduces to solving these problems individually, so that, without loss of generality, we may assume that $m=1$, $A=J_{n}(\lambda)$ with $\re \la >0$ and $B=(b_{r-s})_{r,s=1}^n=\sum_{j=0}^{n-1}b_j J_n(0)^j$ is an upper triangular Toeplitz matrix, which we will do in the remainder of the proof.  

Note that for $f\in\BC(z)$ without a pole at $\la$, we have that 
\[
f(A)=f(J_n(\la))=\sum_{j=0}^{n-1}\frac{f^{(j)}(\la)}{j!} J_n(0)^j.
\]
Hence we need to show that $f$ in \eqref{fz}  satisfies $f^{(j)}(\la)=j! b_j$ for $j=1,\ldots, n-1$. Note that 
\[
v^* e_1=e_1^* e_1=1,\  A e_1= \la e_1,\quad 
v^* e_j=0,\  A e_j=\la e_j + e_{j-1}\ \ \mbox{for}\ \  j>1. 
\]
Therefore, we have that 
\[
(\la H- HA + B^*vv^*)e_1=B^*v \ \ \mbox{and}\ \ 
(\la H- HA + B^*vv^*)e_j = -He_{j-1},\ \ j>1.
\]
We then obtain that 
\[
b_0=v^*B e_1 =v^*B[zH-HA+B^*vv^*]^{(-1)} B^*v=f(\la),
\]
and, recursively, for $j=2,\ldots,n$ and by using Lemma \ref{deriv}, we have that 
\begin{align*}
b_{j-1}&=v^* B e_j 
= -v^* B [\la H-HA+B^*vv^*]^{(-1)} H e_{j-1}\\
&= (-1)^{j-1}v^* B \left( [\la H-HA+B^*vv^*]^{(-1)} H\right)^{j-1} e_1\\
&= (-1)^{j-1}v^* B \left( [\la H-HA+B^*vv^*]^{(-1)} H\right)^{j-1} \times\\
&\qquad\qquad \times[\la H-HA+B^*vv^*]^{(-1)}B^*v\\ 
&= \frac{f^{(j-1)}(\la)}{(j-1)!}.
\end{align*}
The powers of the generalized inverses make sense, since $x\in \kr H \cap \kr \Ga$ implies that $x\perp \im H$. The last identity follows from Lemma \ref{deriv}.   
\end{proof}

When a solution $f\in\cP\cO^{m \times m}$ is given by a realization \eqref{F}, then the matrix $H$ in \eqref{Lya} can be written explicitly in terms of the realization.

\begin{lemma}\label{lem}
Let $A,B\in\BC^{n \times n}$ with $\re A>0$. When $f\in\cP\cO^{m \times m}$ is given via \eqref{F} and satisfies $f(A)=B$, then $H$ in \eqref{Lya} is given by
\begin{equation}\label{H}
H=(C^* \otimes I_n)(I_r \otimes A - i\Lambda \otimes I_n)^{*-1}(I_r \otimes vv^* )(I_r \otimes A - i\Lambda \otimes I_n)^{-1}(C \otimes I_n) + K\otimes vv^*.
\end{equation}
\end{lemma}

\begin{proof}
    If $f$ is given by \eqref{F}, then
    $$ B= f(A) = D\otimes I_n + (C^* \otimes I_n)(I_r \otimes A - i\Lambda \otimes I_n)^{-1}(C \otimes I_n) + K\otimes A.$$
 Using this, it is straightforward to check that $H$ in \eqref{H} satisfies \eqref{Lya}.
\end{proof}

Note that this  lemma provides a different proof of Theorem \ref{ot45}(ii)$\Rightarrow$(iii) as $H$ in \eqref{H} is clearly positive semidefinite.

Let us illustrate Theorem \ref{ot45} with a couple of examples.

\begin{example}\rm
    Let $A=B=\begin{pmatrix}
        1 & 0 \cr 0 & 2
    \end{pmatrix}$, then $v=\begin{pmatrix}
        1  \cr 1
    \end{pmatrix}$.  Solving \eqref{Lya}, we find
    $H=\begin{pmatrix}
        1 & 1 \cr 1 & 1
    \end{pmatrix}$, which is positive semidefinite. The solution \eqref{fz} gives $$f(z)= \begin{pmatrix} 1 & 2 \end{pmatrix} 
    \begin{pmatrix} z & z-1 \cr z+1 & z \end{pmatrix}^{-1} \begin{pmatrix} 1 \cr 2 \end{pmatrix} =z,$$ which is indeed an interpolant in ${\mathcal PO}$.
\end{example}

\begin{example}\label{ex2}\rm
(same as \cite[Section 4, Problem 1]{YS67}).
    Let $A=\begin{pmatrix}
        1 & 0 & 0\cr 0 & 2 &0
\cr 0 & 0 & 3    \end{pmatrix}$ and $B=\begin{pmatrix}
        2 & 0 & 0\cr 0 & \frac52 &0
\cr 0 & 0 & \frac{10}{3}    \end{pmatrix}$. We also have $v=\begin{pmatrix}
        1  \cr 1 \cr 1
    \end{pmatrix}$. Solving \eqref{Lya} gives $H=\begin{pmatrix}
        2 & \frac32 & \frac43 \cr \frac32 & \frac54 & \frac76 \cr \frac43 & \frac76 & \frac{10}{9}
    \end{pmatrix}$, which is positive semidefinite, 
    and $$ HA-Bvv^*=\begin{pmatrix}
        0 & 1 & 2 \cr -1 & 0  & 1  \cr -2 & -1 & 0
    \end{pmatrix}.
    $$ Now
    $$(zH-HA+Bvv^*)^{-1} = $$ $$\frac{1}{4z}\begin{pmatrix} z^2+36 & -4(z^2- z + 18) &  3(z^2 - 2 z + 12) \cr
-4(z^2 + z + 18)  & 16(z^2 + 9) & -12(z^2 - z + 6) \cr
3(z^2 + 2z + 12) & -12(z^2 + z + 6) &  9(z^2 + 4) \end{pmatrix}.$$ 
The solution \eqref{fz} gives $f(z)=z+\frac{1}{z}$, which is indeed an interpolant in ${\mathcal PO}$.
\end{example}

%
%


\subsection{Condition for a unique solution in the case $m=1$}

When $H$ is positive definite then there exist infinitely many solution to the interpolation problem; see, e.g., \cite[Section 22]{BGR}, where a description of the set of all solutions is given. 
In the scalar valued case ($m=1$), the solution is unique as soon as $A$ is nonderogatory and $H$ is singular, as we state next.

\begin{theorem}\label{unique} Let $A\in {\mathbb C}^{n\times n}$ be nonderogatory with $\re A >0$, and let $B\in {\mathbb C}^{n\times n}$. There is a unique $f \in \PO$ such that $f(A)=B$ if and only if the solution $H$ to $$ A^*H+HA = B^*vv^* + vv^*B $$ is singular positive semidefinite. This unique solution is given by \eqref{fz}.
\end{theorem}

When $A$ is not nonderogatory one has automatically redundancy built into the interpolation problem and $H$ will be singular regardless of uniqueness or not. For instance, if $A=B=I_2$ then $H=vv^*$ is singular positive semidefinite, but any convex combination of $z$ and $\frac{1}{z}$ is a solution.

We need the following lemma for the proof of Theorem \ref{unique}.

\begin{lemma}\label{linind} Let $A \in {\mathbb C}^{n\times n}$ be nonderogatory with $\re A >0$ and let $v$ be as in \eqref{v}. Then for different $\omega_j \in {\mathbb R}$, $j=1,\ldots , r$, with $r \le n$,
the vectors $$ (A-i\omega_1)^{*-1}v, \ldots , (A-i\omega_r)^{*-1}v ,$$ are linearly independent. Also, for $r <n$, the vectors
$$v , (A-i\omega_1)^{*-1}v, \ldots , (A-i\omega_r)^{*-1}v ,$$ are linearly independent. 
\end{lemma}

\begin{proof} 
Let $r<n$. We first observe that the matrices 
$$ (A-i\omega_1)^{*-1}, \ldots , (A-i\omega_r)^{*-1}$$
are linearly independent. Indeed, suppose that 
$$ \sum_{j=1}^r \alpha_j (A-i\omega_j)^{*-1} = 0.$$ Then, multiplying with $\prod_{j=1}^r (A-i\omega_j)^*$, we get that
$p(A^*)=0$, where $p(z)$ is a polynomial of degree $r-1<n$. But the minimal polynomial of $A^*$ is of degree $n$, and thus we must have that $p(z)\equiv 0$. This yields that $\alpha_j=0$, $j=1,\ldots , r$. Due to (the proof of) Lemma \ref{vprop}, $ \sum_{j=1}^r \alpha_j (A-i\omega_j)^{*-1} = 0$ if and only if $ \sum_{j=1}^r \alpha_j (A-i\omega_j)^{*-1}v = 0$, and thus the first statement is proven. The proof of the second statement is similar.
\end{proof}

\begin{proof}[Proof of Theorem \ref{unique}]
First, we observe that \cite[Section 22]{BGR} implies that if $H$ is invertible, then there are infinitely many solutions. Thus, it remains to show that the solution is unique when $H$ is singular positive semidefinite. We follow some of the ideas from \cite{YS67}. 

Let $x_1, \ldots , x_k$ be a basis for ${\rm Ker} H$ and put $X=\begin{pmatrix} x_1 & \cdots & x_k \end{pmatrix}.$ If $f \in \PO$ with realization \eqref{F} is a solution, then \eqref{H} gives that $H$ is a sum of positive semidefinite rank one matrices as follows
\begin{equation}\label{H2} H= \alpha vv^* + \sum_{j=1}^r a_j (A-i\omega_j)^{*-1}vv^*(A-i\omega_j)^{-1},\end{equation} where there are $n-k$ non-trivial terms, $\alpha \ge 0$, and $a_j>0$. Due to Lemma \ref{linind} the rank 1's appearing are linearly independent, and thus we have that $\alpha >0$ if and only if $v\in {\rm Ran} H = ({\rm Ker} H)^\perp$, so this happens if and only if $v^*X=0$. Similarly, we need that 
$$ v^* (A-i\omega_j)^{-1} X=0,\quad j=1,\ldots , r.$$
The equations
$$ v^* (A-i \omega)^{-1} x_j=0,\quad \omega \in {\mathbb R},\ j=1,\ldots, k, $$
lead to polynomial equations $p_j(i\omega):= \det(i\om I-A)v^* (A-i \omega)^{-1} x_j=0$, $j=1,\ldots, k$, where $\deg p_j \le n-1$ and $\deg p_j < n-1$ when $v^*X=0$. These equations determine the poles $i\omega_1, \ldots , i\omega_r$ uniquely, with $r \le n-1$ when $v^*X\neq 0$ and $r<n-1$ when $v^*X=0$. Hence $\om_1,\ldots,\om_r$ are uniquely determined by the data.

Finally, one needs to observe that due to Lemma \ref{linind}, the matrices $$  vv^* , (A-i\omega_1)^{*-1}vv^*(A-i\omega_1)^{-1}, \ldots , (A-i\omega_r)^{*-1}vv^*(A-i\omega_r)^{-1} ,$$ or equivalently, the vectors $$v, (A-i\omega_1)^{*-1}v, \ldots , (A-i\omega_r)^{*-1}v ,$$ are linearly independent (where $v$ only appears when $v^*X=0$).



The linear independence shows that the decomposition \eqref{H2} is unique, that is, the numbers $a_j$ are uniquely determined from $H$, and hence from $A$ and $B$. Since we already observed that $\om_1,\ldots,\om_r$ are uniquely determined by $A$ and $B$, if follows that the right hand side of \eqref{F} is uniquely determined by $A$ and $B$, and hence the solution $f$ is unique.
\end{proof}

The uniqueness proof suggests an alternative way to construct a solution.
If $H$ is singular, let $0\neq x \in {\rm Ker} H$ and let
$$ g(z)=v^* (zI-A)^{-1} x =: \frac{p(z)}{\det (zI-A)}.$$
Note that $\deg p(z) \le n-1$ and $\deg p(z) < n-1$ when $v^*x=0$. Let the roots of $p(z)$ on $i{\mathbb R}$ be denoted by $i\omega_j$, $j=1,\ldots, r$, where we include $i\infty$ when $v^*x=0$. Now write $B$ as a linear combination of
\[
 I,\  (A-i \omega_1I_n)^{-1}, \ldots , (A-i \omega_rI_n)^{-1} ,\]
 where $A:=(A-i \infty I_n)^{-1}$. That $B$ can be written in this way follows from the formula of $B=f(A)$ in the proof of Lemma \ref{lem}.
If $B=i a_0I + \sum_{j=1}^r {a_j} (A-i \omega_jI_n)^{-1} $ then the unique solution is
$$ f(z)= i a_0  +\sum_{j=1}^r  a_j (z-i \om_j)^{-1},$$
where $(z-i \infty I_n)^{-1}$ should be interpreted as $z$. 
The above results imply that $a_j\ge 0$, $j=1,\ldots r.$

When $H$ is invertible, one can choose a pole $i\om$ on $i{\mathbb R}$ and adjust the interpolation data such that the solution includes a positive multiple of $(z-i \om_j)^{-1}$. As long as the new $H (=H_{old}-a(A-i\omega)^{*-1}vv^*(A-i\omega)^{-1})$ remains positive definite, one can repeat this, and when it becomes singular positive semidefinite one continues as in the previous paragraph.

Let us apply the above to Example \ref{ex2}.

\noindent{\bf Example \ref{ex2} continued.} We had $H=\begin{pmatrix}
        2 & \frac32 & \frac43 \cr \frac32 & \frac54 & \frac76 \cr \frac43 & \frac76 & \frac{10}{9}
    \end{pmatrix}$, and $v=\begin{pmatrix}
        1\cr 1\cr 1
    \end{pmatrix}$. We have that 
    $ {\rm Ker} H = {\rm span} \{ x \}$, where $x=\begin{pmatrix}
        1 \cr -4 \cr 3
    \end{pmatrix}.$
We see that $v^*x=0$, so the solution will have a pole at infinity. We see that the maximal $\alpha$ such that $H-\alpha vv^* \ge 0$ is $\alpha=1$, yielding a pole at infinity of the form $\alpha z=z$. Next, we consider the expression
$$ v^* (A-z)^{-1} x = \frac{1}{1-z} -\frac{4}{2-z}+\frac{3}{3-z} = $$ $$\frac{(2-z)(3-z)-4(1-z)(3-z)+3(1-z)(2-z)}{(1-z)(2-z)(3-z)} = \frac{2z}{(1-z)(2-z)(3-z)}.$$
The only root in $i{\mathbb R}$ of the numerator is 0. Thus, the solution will have a pole at 0. It remains to determine the maximal $\beta$ such that
$$ H-\beta (A^*-0)^{-1}vv^*(A-0)^{-1} \ge 0,$$ and we find the maximal to be $\beta =1$. Thus the pole at 0 is of the form $\beta \frac{1}{z}= \frac{1}{z}$. Thus we find that
$f(z)=z+\frac{1}{z}$ is the unique solution.


\subsection{Discussion on the description of all solutions}

To end this section, let us comment on descriptions of the set of all solutions. When $H$ is positive definite, then \cite[Theorem 22.2.2]{BGR} provides an explicit construction of a linear fractional map
\begin{equation}\label{f} f = (\theta_{11} + \theta_{12}g)(\theta_{21} + \theta_{22} g)^{-1},\end{equation}
that is a one-to-one map from $g\in\cP^{m\times m}$ without poles on $i\BR$ onto the set of solutions $f\in\cP^{m\times m}$ without poles on $i\BR$. One may rewrite \eqref{f} as 
\begin{equation}\label{f2} f=L_{11}+ L_{12} g (I + L_{22} g )^{-1} L_{21}, \end{equation}
where
\begin{align*}
 L_{11} &=
(I_n\otimes v^*)B[zH-(H(I_n\otimes A)-B^*(I_n\otimes vv^*))]^{(-1)} B^*(I_n\otimes v),\\
L_{12} &= I+ (I_n\otimes v^*)B[zH-(H(I_n\otimes A)-B^*(I_n\otimes vv^*))]^{(-1)} (I_n\otimes v),\\
L_{21}&= I-(I_n\otimes v^*)[zH-(H(I_n\otimes A)-B^*(I_n\otimes vv^*))]^{(-1)} B^*(I_n\otimes v),\\
 L_{22} &= (I_n\otimes v^*)[zH-(H(I_n\otimes A)-B^*(I_n\otimes vv^*))]^{(-1)} (I_n\otimes v).
\end{align*}
We expect that the adjusted map \eqref{f2} can be extended to $g\in\cP^{m\times m}$, mapping $g\in\cP\cO^{m\times m}$ to $f\in\cP\cO^{m\times m}$, giving a description of all solutions, also when $H$ is singular. Such a result would require $g$ to range in a subset of $\cP^{m\times m}$, such that \eqref{f2} is well-defined. Technical details would need to be worked out and are expected to be quite intricate, as similar technicalities appear in the matrix valued Nevanlinna-Pick problem when the Pick matrix is singular.

\newpage

\section{The one variable case and completely positive maps}

We first present a `coordinate free' version of Theorem \ref{ot45}, which provides further details on when there are solutions with a pole at $\infty$. Let $w$ be the vector defined by
\begin{equation}\label{w}
w:=\sum_{i=1}^n e_i^{(n)} \otimes e_i^{(n)} \in\BC^{n^2},
\end{equation}
and note that we have the rank-1 factorization $\left[E_{ij}\right]_{i,j=1}^m=ww^*$, where $E_{ij} \in {\mathbb C}^{n\times n}$ has a 1 in position $(i,j)$ and zeros elsewhere. In this section we will set $e_i = e_i^{(n)}$; in other words, for convenience we will leave off the superscript $(n)$. 

Before stating our first result, we need to introduce some notation and terminology. Let ${\rm PSD}_k$ denote the cone of $k\times k$ positive semidefinite matrices. We call a linear map $\Phi : \BC^{n \times n} \to \BC^{l \times l}$ {\it positive} if $\Phi ({\rm PSD}_n) \subseteq {\rm PSD}_l$.  Let $\textup{id}_{\BC^{k \times k}} \otimes \Phi$ be the linear map induced by $\Phi$ in the following way
\[
\pmat{A_{11} & \hdots & A_{1k} \\ \vdots & \ddots & \vdots \\ A_{k1} & \hdots & A_{kk}} \mapsto \pmat{\Phi\left(A_{11}\right) & \hdots & \Phi\left(A_{1k}\right) \\ \vdots & \ddots & \vdots \\ \Phi\left(A_{k1}\right) & \hdots & \Phi\left(A_{kk}\right)}.
\] 
We call $\Phi$ {\em $k$-positive} in case $\textup{id}_{\BC^{k \times k}} \otimes \Phi$ is positive, and {\em completely positive} if $\Phi$ is $k$-positive for each $k\in \BN$. Complete positivity is characterized by Choi's Theorem.

\begin{theorem}[\cite{Choi}]\label{T:Choi}
For a linear map $\Phi : {\mathbb C}^{n \times n} \to {\mathbb C}^{l \times l}$, the following are equivalent:
\begin{itemize}
\item[(i)] $\Phi$ is completely positive. 

\item[(ii)] The {\em Choi matrix} $$\BP:=\left[\Phi(E_{ij})\right]_{i,j=1}^n=(\textup{id}_{\BC^{n \times n}} \otimes \Phi) (ww^*)$$ of $\Phi$ is in ${\rm PSD}_{nl}$.

\item[(iii)] $\Phi$ admits a {\em Kraus representation}:\ There exist matrices $S_1,\ldots,S_s\in\BC^{n \times l}$ such that for all $X\in\BC^{n \times n}$
\begin{equation} \label{D:Kraus}
\Phi(X)=\sum_{r=1}^s S_r^* X S_r =\bS^* (I_s \otimes X) \bS,\ \ \mbox{with }\bS=\sbm{S_1\\ \vdots\\ S_s}.
\end{equation}
\end{itemize} 
\end{theorem}


We now arrive at the main result of this section. We use that, by Lyapunov's Theorem, we have that $\{ XA+A^*X : X \in {\mathbb C}^{n\times n} \} = {\mathbb C}^{n\times n}$ when $\re A >0$.

\begin{theorem} \label{T:Main1}
Let $A\in\BC^{n \times n}$ with $\re A > 0$ and $ B\in\BC^{nm \times nm}$. Then the following statements are equivalent: 
\begin{itemize}
\item[(i)] $f(A)=B$ for a function $f \in \PO^{m \times m}$.

\item[(ii)] $f(A)=B$ for a function $f \in \cP^{m \times m}$.

\item[(iii)] $\cL: \BC^{n \times n} \to \BC^{nm \times nm}$ defined via
$$ \cL (XA+A^*X):=(I_m\otimes X) B + B^*(I_m \otimes X), \ X\in \BC^{n \times n},$$ 
is a completely positive linear map.

\item[(iv)] The Choi matrix 
\[
\BP_{A,B}:=\left[\cL\left(E_{ij}\right)\right]_{i,j=1}^n=(\textup{id}_{\BC^{n \times n}} \otimes {\mathcal L})(ww^*)
\]
is positive semidefinite.

\end{itemize}
Assume (i)-(iv) hold. Then the kernel inclusion condition
\begin{equation}\label{KerIncl} 
 \kr\left[ I_m\otimes E_{ij}+\cL\left(E_{ij}\right)\right]_{i,j=1}^n \subseteq \kr\left[  B^*(I_m \otimes E_{ij})B + \cL\left(A^*E_{ij} A\right)\right]_{ij=1}^n , \end{equation}
holds if and only if there is a solution of the form
\begin{equation}\label{SolRealProper} 
f(z)=D+ \ell^*\left(zI_r - M\right)^{-1}\ell, \mbox{ with } D=-D^* ,\ M^*=-M,
\end{equation}
and in that case a solution $f$ as in \eqref{SolRealProper} exists with $r=\rank \BP_{A,B}$. Furthermore, all solutions are of the form \eqref{SolRealProper} if and only if 
\begin{equation}\label{ranran}\cM:= \ran W^*
\cap \ran \BP_{A,B} = \{ 0 \},\end{equation} where $W=\sum_{i=1}^n e_i^* \otimes I_m \otimes e_i^*$ If the kernel inclusion \eqref{KerIncl}  does not hold then all solutions are of the form 
\begin{equation}\label{SolRealNonproper}
\begin{aligned}
f(z)&=K z+D +\ell^*\left(zI_r-M\right)^{-1}\ell,\\ 
&\mbox{ with } K\ge 0, K\neq 0,\ D=-D^*,\ M^*=-M,
\end{aligned}
\end{equation}
and $K\le K_{\rm max}$, where $W^*K_{\rm max}W$ is the Schur complement of $\BP_{A,B}$ supported on $\cM$. Moreover, a solution with $K=K_{\rm max}$ and $r=\rank \BP_{A,B} - \rank K_{\max}$ exists.
Finally, if in addition $A$ and $B$ are real, then the matrices $K$, $D$, $\ell$ and $M$ can be chosen to be real (and thus $f\in \PRO^{m\times m}$).
\end{theorem}

The following lemmas will be crucial in the proof of Theorem \ref{T:Main1}(iii)$\to$(i). Let us denote 
$$ W= \sum_{i=1}^n e_i^* \otimes I_m \otimes e_i^* = \begin{pmatrix} I_m\otimes e_1^* & \cdots & I_m\otimes e_n^* \end{pmatrix}. $$

\begin{lemma}\label{L:KNLequations}
Let $A,B\in\BC^{n \times n}$ with $\re A > 0$ and assume that $\cL$ in Theorem \ref{T:Main1} is completely positive with Kraus representation $\cL (Y)= \sum_{i=1}^r C_i^*YC_i$ and Choi matrix $\BP_{A,B}$. Define $K\in \BC^{r \times mn^2}$ and $L,N \in \BC^{(m + r)\times mn^2}$ by
$$  K=\begin{pmatrix}
     w^*(I_n\otimes C_1) \cr \vdots \cr w^*(I_n\otimes C_r)\end{pmatrix},  L=\begin{pmatrix} W
    \cr w^*(I_n\otimes C_1) \cr \vdots \cr w^*(I_n\otimes C_r) 
    
\end{pmatrix}, $$ $$ N= \begin{pmatrix} W(I_n \otimes B)
    \cr -w^*(I_n\otimes AC_1) \cr \vdots \cr -w^*(I_n\otimes AC_r)
    
\end{pmatrix},$$
 with $w \in \BC^{n^2}$ as in \eqref{w}. Then we have
\begin{equation}\label{KNLmatEQS}
\BP_{A,B}=K^*K,\quad N^*L=-L^*N. 
\end{equation}
Moreover, $K$ is full row-rank if and only if the Kraus representation $\cL (Y)= \sum_{i=1}^r C_i^*YC_i$ is minimal. Furthermore, \eqref{KerIncl} holds if and only if $\kr L \subseteq \kr N$, and in that case there exists an $S=-S^* \in {\mathbb C}^{(m+r)\times (m+r)}$ such that $SL=N$.
\end{lemma}

\begin{proof}

Since $\rank \BP_{A,B}=r$, the factorization \begin{align*} 
 \BP_{A,B}&=\left[ \cL(E_{ij} )\right]_{i,j=1}^n = \left[\sum_{k=1}^r C_k^* E_{ij} C_k \right]_{i,j=1}^n\\&= \sum_{k=1}^r (I_n \otimes C_k^*) \left[E_{ij}\right]_{i,j=1}^n (I_n\otimes C_k)\\&= \sum_{k=1}^r (I_n \otimes C_k^*) ww^*(I_n\otimes C_k)= K^*K,
\end{align*}
is a minimal rank factorization of $\BP_{A,B}$, i.e., $K$ has full row-rank, exactly when the Kraus representation $\cL (Y)= \sum_{i=1}^r C_i^*YC_i$ is minimal.

Next, using that for all $i,j=1,\ldots , n$ we have 
 $$\cL\left(E_{ij}A + A^*E_{ij} \right) = (I_m \otimes E_{ij})B + B^*(I_m \otimes E_{ij}),$$ we observe that \begin{align*}
    0&= \begin{pmatrix}
         (I_m \otimes E_{ij})B+B^*(I_m \otimes E_{ij}) - \cL\left(E_{ij}A+A^*E_{ij}\right)
     \end{pmatrix}_{i,j=1}^n \\&= \begin{pmatrix}
         (I_m \otimes E_{ij})B+B^*(I_m \otimes E_{ij}) - \cL\left(E_{ij}A\right) -\cL\left(A^*E_{ij}\right)
     \end{pmatrix}_{i,j=1}^n
 \end{align*} which results in \begin{align*}
     \begin{pmatrix}
         B^*(I_m \otimes E_{ij}) - \cL\left(A^*E_{ij}\right)
     \end{pmatrix}_{i,j=1}^n=-\begin{pmatrix}
         (I_m \otimes E_{ij})B -\cL\left(E_{ij}A\right) 
     \end{pmatrix}_{i,j=1}^n.
 \end{align*} From this we obtain the matrix equation
\begin{align} 
N^* L &= (I_n \otimes B^*) [(I_m\otimes E_{ij})_{i,j=1}^n] - \sum_{k=1}^r (I_n \otimes C_k^* A^* ) ww^* (I_n \otimes C_k) \notag \\ 
&= \begin{pmatrix}
         B^*(I_m \otimes E_{ij}) -  \sum_{k=1}^r  C_k^* A^* E_{ij} C_k
     \end{pmatrix}_{i,j=1}^n = \begin{pmatrix}
         B^*(I_m \otimes E_{ij}) - \cL\left(A^*E_{ij}\right)
     \end{pmatrix}_{i,j=1}^n\notag \\
&= -\begin{pmatrix}
         (I_m \otimes E_{ij})B -\cL\left(E_{ij}A\right) 
     \end{pmatrix}_{i,j=1}^n  
=-L^*N.\label{N*L workings} 
\end{align} 
Hence the identities in \eqref{KNLmatEQS} hold. 

Next observe that
\begin{align*}
     L^*L&= [(I_m\otimes E_{ij})_{i,j=1}^n] + K^*K=\left[ I_m\otimes E_{ij}+ \cL\left(E_{ij}\right)\right]_{i,j=1}^n, \quad \end{align*} 
and 
\begin{align} 
N^*N
&= (I_n \otimes B^*)[(I_m\otimes E_{ij})_{i,j=1}^n](I_n \otimes B) + \sum_{k=1}^r (I_n \otimes C_k^*A^*)ww^*(I_n \otimes AC_k) 
      \notag \\
&= \begin{pmatrix}
         B^*(I_m \otimes E_{ij})B
     \end{pmatrix}_{i,j=1}^n + \begin{pmatrix}
         \sum_{k=1}^r C_k^*A^*E_{ij}AC_k
     \end{pmatrix}_{i,j=1}^n \notag \\
&=\left[  B^*(I_m \otimes E_{ij})B + \cL\left(A^*E_{ij} A\right)\right]_{ij=1}^n.\notag
\end{align}
It follows that \eqref{KerIncl} can be rewritten as $\kr L^*L \subseteq \kr N^*N$. In general we have that $\kr G^*G = \kr G$ for any matrix $G$. Using this, we see that \eqref{KerIncl} is equivalent to $\kr L \subseteq \kr N.$ 

Finally, assume that the matrices $L$ and $N$ in Lemma \ref{L:KNLequations} satisfy $\kr L \subseteq \kr N.$ By Douglas' Lemma there exists a unique matrix $S \in {\mathbb C}^{(m+r)\times (m+r)}$ such that 
 $$ S L=N \quad \mbox{and}\quad (\ran L)^\perp \subseteq \kr S.$$
Replacing $N$ by $SL$ in the identity $N^*L=-L^*N$, it follows that $L^*S^*L=-L^*SL$. Write 
$$ S=\begin{pmatrix} S_{11} & S_{12} \cr S_{21} & S_{22} \end{pmatrix} : \begin{matrix} \ran L \cr \oplus \cr (\ran L)^\perp \end{matrix} \to \begin{matrix} \ran L \cr \oplus \cr (\ran L)^\perp \end{matrix}. $$
Then $L^*S^*L=-L^*SL$ implies that $S_{11}$ is skew-Hermitian. Moreover, $(\ran L)^\perp \subseteq \kr S$ gives that $S_{12}=0$ and $S_{22}=0$.  Now put
$$ \hat S = \begin{pmatrix} S_{11} & -S_{21}^* \cr S_{21} & 0 \end{pmatrix}. $$ Then we have that $\hat S $ is skew-Hermitian and $\hat  S L=SL=N$ (as $S$ and $\hat S$ agree on $\ran L$).
\end{proof}

\begin{lemma}\label{L:winRanP}
Let $A,B,\cL$ and $\BP_{A,B}$ as in Theorem \ref{T:Main1}. If (iii) in Theorem \ref{T:Main1} and \eqref{ranran} hold, then \eqref{KerIncl} holds. 
\end{lemma}

\begin{proof} Recall that
$$ W= \sum_{i=1}^n e_i^* \otimes I_m \otimes e_i^* . $$ Note that
$$W^*W = (I_m\otimes E_{ij})_{i,j=1}^n. $$ Assume (iii) in Theorem \ref{T:Main1} and $ {\rm Ran} W^* \cap \ran \BP_{A,B} = (0)$. Take any minimal Kraus representation  $\cL (Y)= \sum_{i=1}^r C_i^*YC_i$ of $\cL$ and define $K$, $L$ and $N$ as in Lemma \ref{L:KNLequations}. By Lemma \ref{L:KNLequations}, in order to show that \eqref{KerIncl} holds it suffices to show that $\kr L \subset \kr N$. Let $u\in \kr L$. Using \eqref{KNLmatEQS} and \eqref{N*L workings} we obtain that
$$ 0=N^*Lu=-L^*Nu = -W^*W(I_n\otimes B)u - 
K^* \begin{pmatrix} 
      -w^* (I_n\otimes A C_1) \cr \vdots \cr  -w^* (I_n\otimes A C_r)
 \end{pmatrix} u. $$ 
If $W(I_n\otimes B)u\neq 0$, then this identity implies that $(0) \neq \ran W^*W \cap \ran K^* = \ran W^* \cap \ran K^*K = \ran W^* \cap\ran \BP_{A,B}$, contradicting our assumption. Hence $W(I_n\otimes B)u=0$.

Again considering the above identity, $W(I_n\otimes B)u=0$ implies that 
$$ K^* \begin{pmatrix} 
     -w^* (I_n\otimes A C_1) \cr \vdots \cr  -w^* (I_n\otimes A C_r)
 \end{pmatrix} u=0.$$
Since we are working with a minimal Kraus representation,   $K^*$ is injective by Lemma \ref{L:KNLequations}, and hence 
$$ \begin{pmatrix} 
     -w^* (I_n\otimes A C_1) \cr \vdots \cr  -w^* (I_n\otimes A C_r)
 \end{pmatrix} u=0.  $$ 
 Combining this with $W(I_n\otimes B)u= 0$, gives that $Nu=0$, which finishes the proof.
\end{proof}

\begin{lemma}\label{PW}
    Let $P\in{\mathbb C}^{s\times s}$ be positive semidefinite, and $X\in {\mathbb C}^{s\times l}$ with $\rank X = l$. Then there exists a unique $l\times l$ positive semidefinite matrix $K_{\rm max}$ such that
    $$P-XK_{\rm max} X^* \ge 0 ,$$
    and $K_{\rm max}$ is maximal (in the Loewner order) with these properties. In this case $\ran(P-XK_{\rm max} X^*) \cap \ran X =(0).$
\end{lemma}
\begin{proof}
    Let $\cM = \ran P \cap \ran X$. 
    Let $S$ be the Schur complement of $P$ supported on $\cM$ (see \cite{Ando}); i.e., $\ran S \subseteq \cM$, $S\ge 0$ and $P-S\ge 0$ and $S$ is maximal (in the Loewner order) with these properties. Since $\ran S \subseteq \ran X$, $X$ has full column rank and $S\ge 0$, we may write $S$ uniquely as $S=XK_{\rm max} X^*$, for some positive semidefinite matrix $K_{\rm max}$. The matrix $K_{\rm max}$ has the desired properties.
\end{proof}

\noindent {\it Proof of Theorem \ref{T:Main1}.}
(i) $\to$ (ii) is trivial. The equivalence of (iii) and (iv) is Choi's Theorem \cite[Theorem 2]{Choi}.
It remains to show (ii)$\to$[(iii)\&(iv)]$\to$(i).

Assume (ii) holds and let $k\in{\mathbb N}$ and $T=(T_{ij})_{i,j=1}^k \in{\mathbb C}^{nk\times nk}$ be invertible. Then whenever $T(\oplus A)T^{-1}$ has positive definite real part (with $k$ copies of $A$ appearing), it follows that \[ (I_m \otimes T_{ij})_{i,j=1}^k (\oplus B)(I_m \otimes (T^{-1})_{ij})_{i,j=1}^k=\]\[ (I_m \otimes T_{ij})_{i,j=1}^k (\oplus f(A))(I_m \otimes (T^{-1})_{ij})_{i,j=1}^k\] has positive semidefinite real part. Let $P=T^*T$, and write $P=(P_{ij})_{i,j=1}^k$. From the above, we deduce that
\[P (\oplus A) + (\oplus A)^*P >0 \quad {\implies} \quad \] \[(I_m \otimes P_{ij})_{i,j=1}^k (\oplus B) + (\oplus B)^* (I_m \otimes P_{ij})_{i,j=1}^k  \ge 0.\] 
Thus, it follows that if $P>0$ and $P (\oplus A) + (\oplus A)^*P >0$, then $(\textup{id}_{\BC^{k \times k}} \otimes \cL_{A,B})(P (\oplus A) + (\oplus A)^*P) \ge 0.$ But then the same holds true when $P \ge 0$ and $P (\oplus A) + (\oplus A)^*P \ge 0$, by a $P+\epsilon I$ and continuity argument. Next, let $X>0$. Then by Lyapunov's Theorem there exists $P\ge 0,$ such that $X=P (\oplus A) + (\oplus A)^*P.$ But then $(\textup{id}_{\BC^{k \times k}}\otimes \mathcal{L}_{A,B}) (X) \ge 0$. Again, by continuity, we find that $(\textup{id}_{\BC^{k \times k}} \otimes \mathcal{L}_{A,B})(X)\ge 0 $ whenever $X\ge 0$. Thus ${\mathcal L}_{A,B}$ is completely positive and (iii) \& (iv) follow. 

Next, assume that (iii) \& (iv) hold. Since $\cL$ is completely positive, we may 
choose a minimal Kraus representation for it: ${\mathcal L} (Y) =\sum_{j=1}^r C_j^* Y C_j,  $ with $C_j$ of size $n\times nm$.
First we assume that \eqref{KerIncl} holds.
Then we can apply Lemma \ref{L:KNLequations}, where for
$$ L=\begin{pmatrix}
    W \cr w^*(I_n\otimes C_1) \cr \vdots \cr w^*(I_n\otimes C_r)
    
\end{pmatrix}, N= \begin{pmatrix} W(I_n\otimes B)
    \cr -w^*(I_n\otimes AC_1) \cr \vdots \cr -w^*(I_n\otimes AC_r)
    
\end{pmatrix},$$
we find a $(m+r)\times (m+r) $ matrix 
$$ S=-S^* =: \begin{pmatrix}
    D & \ell^* \cr -\ell & -M 
\end{pmatrix},$$
such that $SL=N$. Here $D=-D^* \in {\mathbb C}^{m\times m} $ and $M=-M^* \in {\mathbb C}^{r\times r} $. Define now $f(z)$ via \eqref{SolRealProper}. By Theorem \ref{PO}, $f(z)\in \PO^{m\times m}$. 

It remains to check that $f(A)=B.$
For this, observe that the first row of the equality $S L=N$ gives \[D\begin{pmatrix} I_m\otimes e_1^* & \cdots & I_m\otimes e_n^* \end{pmatrix} + \ell^*K=\begin{pmatrix}
     (I_m\otimes e_1^*)B & \cdots & (I_m\otimes e_n^*)B \end{pmatrix},\] where $K$ is as in Lemma \ref{L:KNLequations}. Writing $\ell^*=\begin{pmatrix} {\ell}_1^* & \hdots & {\ell}_r^*\end{pmatrix}$ it then follows that $$ D (I_m\otimes e_j^*) + \sum_{k=1}^r {\ell}_k^* e_j^* C_k = (I_m\otimes e_j^*) B, \quad \text{for all \,}  j=1,\ldots , n. $$
Therefore
\begin{equation}\label{expressionB} D\otimes I_n + (\ell^*\otimes I_n)\begin{pmatrix} C_1\cr \vdots \cr C_r \end{pmatrix}= D\otimes I_n + \sum_{k=1}^m (\ell^*)_k C_k = B. \end{equation} For $M=[m_{ij}]$ the next rows of the equality $ S L=N$ give that
$$ -\ell_s (I_m \otimes e_j^*) - \sum_{k=1}^r m_{sk} e_j^* C_k = -e_j^* A C_s, \quad \text{for all\, } j=1,\ldots , n,$$ and thus 
$$ -\ell_s \otimes I_n - \sum_{k=1}^r m_{sk} C_k = -A C_s \quad \text{for all\, } s=1,\ldots , r.$$ 
Rewriting the latter equation gives
$$ -\ell\otimes I_n - (M\otimes I_n) \begin{pmatrix} C_1\cr \vdots \cr C_r \end{pmatrix} = -(I_r \otimes A) \begin{pmatrix} C_1\cr \vdots \cr C_r \end{pmatrix}, $$ which combined with \eqref{expressionB} yields
$$ D \otimes I_n + (\ell^*\otimes I_n) \left((I_r \otimes A) - (M\otimes I_n)\right)^{-1} (\ell\otimes I_n)=B.$$
This shows that $f$ in \eqref{SolRealProper} satisfies $f(A)=B$. Thus, $f$ solves the interpolation problem. 

Next, assume that \eqref{KerIncl} does not hold. By Lemma \ref{PW} with $P={\mathbb P}_{A,B}$ and $W=\sum_{i=1}^n e_i^* \otimes I_m \otimes e_i^*$ we obtain a matrix $K_{\rm max}$. Let us now consider the interpolation problem where we want $f(A)=B-K_{\rm max}\otimes A$. In this setting, we have that 
$$ {\mathbb P}_{A, B-K_{\rm max}\otimes A} = {\mathbb P}_{A, B} - W^*K_{\rm max}W,$$
and due to Lemma \ref{PW} (with $X=W^*$) we have that for the problem where we want $f(A)=B-K_{\rm max}\otimes A$, equation \eqref{ranran} holds. Thus, we have that \eqref{KerIncl} is satisfied. Consequently,  we may apply the previous paragraph and obtain a solution of the form \eqref{SolRealProper} such that $f(A)=B-K_{\rm max}\otimes A$. But then we find a solution for the original interpolation problem of the form \eqref{SolRealNonproper} with $K=K_{\rm max}$. Thus we have proven [(iii) \& (iv)]$\to$(i).

Finally, it follows from the argument above, that if \eqref{SolRealNonproper} is a solution, then ${\mathbb P}_{A, B} - W^*KW \ge 0$. But then Lemma \ref{PW} gives that $0\le K \le K_{\rm max}$. In addition, we have that $K_{\rm max}=0$ if and only if \eqref{ranran} holds. This yields the statement that all solutions are of the form \eqref{SolRealProper} if and only if \eqref{ranran} holds.
\hfill $\qed$

\section{Single point interpolation on noncommutative half-planes}
\label{sec:Q-interpolation}

We begin with some background on  noncommutative (abbreviated as 'nc') rational matrix
functions. More details can be found in the monograph \cite{KV}.

A nc rational matrix expression $R$ over a field
$\mathbb K$ is any expression obtained from noncommuting indeterminates
$z_1,\ldots,z_d$, a constant $1\in\mathbb K$, and matrices of compatible
sizes by successive addition, multiplication, and inversion. Thus, for
example, an $\alpha\times\beta$ nc matrix polynomial
has the form
\[
P=\sum_{w\in\mathcal G_d:\,|w|\leq L}P_w z^w,
\]
where $\mathcal G_d$ is the free monoid on the generators
$g_1,\ldots,g_d$ and with identity element denoted by $\emptyset$, the coefficients $P_w$ are $\alpha\times\beta$ matrices,
and, for a word $w=g_{i_1}\cdots g_{i_N}$, we write
$z^w=z_{i_1}\cdots z_{i_N}$ and $|w|=N$, with
$z^\emptyset=1$.

For a $d$-tuple $Z=(Z_1,\ldots,Z_d)$ of $n\times n$ matrices one evaluates
\[
P(Z)=\sum_{w\in\mathcal G_d}P_w\otimes Z^w,
\qquad
Z^w=Z_{i_1}\cdots Z_{i_N},
\]
with $Z^\emptyset=I_n$.  An nc rational matrix expression $R$ can be evaluated at
$Z$ whenever every formal matrix inverse occurring in the expression can
be replaced by an inverse of an actual matrix at $Z$.  This gives the
level-$n$ domain $\operatorname{dom}_n R$, and
\[
\operatorname{dom}R=\coprod_{n=1}^\infty\operatorname{dom}_nR.
\]

Two nc rational matrix expressions are called
equivalent if they agree on the intersection of their domains and that
intersection is nonempty.  An equivalence class is called a
nc rational matrix function.  We write $f\in\mathfrak R$
if the nc rational matrix expression $f$ represents the rational function $\mathfrak R$, and
set
\[
\operatorname{dom}_n\mathfrak R
=
\bigcup_{R\in\mathfrak R}\operatorname{dom}_nR,
\qquad
\operatorname{dom}\mathfrak R
=\coprod_{n=1}^\infty\operatorname{dom}_n\mathfrak R.
\]

An nc rational matrix function $f$ is said to be \emph{regular at $0$} if
$0\in\operatorname{dom}_1 f$; equivalently, it admits a convergent
noncommutative power-series expansion in a neighborhood of $0$,
\[
f(z)=\sum_{w\in\mathcal G_d}f_w z^w.
\]
For an nc rational matrix function regular at $0$, evaluation of this series on
$n\times n$ matrices converges in a sufficiently small neighborhood of
$0$ for each fixed $n$.

We now introduce the general noncommutative half-plane associated with a
self-adjoint linear pencil. Fix a positive integer $r$ and matrices
\[
L_1,\ldots,L_d\in\mathbb C^{r\times r},
\qquad L_k=L_k^*.
\]
For $Q(z)=\sum_{k=1}^d z_kL_k$, we define as usual
\begin{equation}\label{eq:Qdef}
Q(H)=\sum_{k=1}^d L_k\otimes H_k,
\qquad
H=(H_1,\ldots,H_d)\in(\mathbb C^{n\times n})^d,
\end{equation}
and, for $n\geq1$,
\begin{equation}\label{eq:HQn}
\mathbb H_{Q,n}
=
\left\{
H\in(\mathbb C^{n\times n})^d:
\operatorname{Re}Q(H)>0
\right\}.
\end{equation}
Here
\[
\operatorname{Re}Q(H)=\frac12\bigl(Q(H)+Q(H)^*\bigr).
\]
We write
\[
\mathbb H_Q=\bigcup_{n\geq1}\mathbb H_{Q,n}
\]
and, since the coefficients $L_k$ are self-adjoint,
\begin{equation}\label{eq:Qstar}
Q(H)^*=Q(H^*),
\qquad
H^*=(H_1^*,\ldots,H_d^*).
\end{equation}


\begin{example}
 {\em The Siegel right half-plane:}
For symmetric matrices
\[
\Sym_n({\mathbb C})
=
\{Z\in{\mathbb C}^{n\times n}:Z^T=Z\},
\]
consider the condition
\[
\Ree(Z)>0.
\]
Writing
\[
Z
=
\sum_{1\leq i\leq j\leq n}
z_{ij}L_{ij},
\]
where
\[
L_{ii}=E_{ii},
\qquad
L_{ij}=E_{ij}+E_{ji}\quad(i<j),
\]
we introduce the Hermitian linear pencil
\[
Q(z)
=
\sum_{1\leq i\leq j\leq n}z_{ij}L_{ij}.
\]
Then the $n\times n$ Siegel right half plane is
\[
\mathbb H_Q
=
\left\{
z\in{\mathbb C}^{n(n+1)/2}:
\Ree Q(z)>0
\right\}.
\]
\end{example}
In addition, one can take direct sums of the above domain (with potentially different sizes). For instance, when we take a direct sum of scalars we get the standard poly-half-plane, which we discuss in further detail in the next section.

\bigskip

Let $\mathcal P^{m\times m}_{\rm nc}(\mathbb H_Q)$ denote the class of
$m\times m$ nc rational matrix functions $f$ whose domain contains
$\mathbb H_Q$ and which satisfy
\begin{equation}\label{eq:positive-real}
\operatorname{Re}f(H)\geq0,
\qquad H\in\mathbb H_Q.
\end{equation}
We write $\mathcal{PO}^{m\times m}_{\rm nc}(\mathbb H_Q)$ for the subclass
satisfying, in addition,
\begin{equation}\label{eq:odd}
f(-H)^*=-f(H^*)
\end{equation}
whenever both sides are defined.

Given
\[
A=(A_{1},\ldots,A_{d})\in\mathbb H_{Q,n},
\qquad
B\in\mathbb C^{nm\times nm},
\]
we consider the interpolation problem
\[
f(A)=B,
\qquad
f\in\mathcal{PO}^{m\times m}_{\rm nc}(\mathbb H_Q).
\]
Put
\begin{equation}\label{eq:Q0}
Q_0:=Q(A)=\sum_{k=1}^dL_k\otimes A_{k}
\in\mathbb C^{rn\times rn}.
\end{equation}
Define the self-adjoint linear space
\begin{equation}\label{eq:SigmaQ}
\Sigma_Q
=
\left\{
(I_r\otimes X)Q_0+Q_0^*(I_r\otimes X):
X\in\mathbb C^{n\times n}
\right\}
\end{equation}
and the linear map $\Phi_Q:\Sigma_Q\to\mathbb C^{nm\times nm}$ by
\begin{equation}\label{eq:PhiQ}
\Phi_Q\left((I_r\otimes X)Q_0+Q_0^*(I_r\otimes X)\right)
=
(I_m\otimes X)B+B^*(I_m\otimes X).
\end{equation}

\begin{theorem}\label{thm:Q-interpolation}
Let $A\in\mathbb H_{Q,n}$ and $B\in\mathbb C^{nm\times nm}$.  If there
exists $f\in\mathcal P^{m\times m}_{\rm nc}(\mathbb H_Q)$ with
$f(A)=B$, then $\Phi_Q$ is completely positive.

Moreover, there exists an interpolant
$f\in\mathcal{PO}^{m\times m}_{\rm nc}(\mathbb H_Q)$ of the form
\begin{equation}\label{eq:F-realization}
f(z)=
D
+\ell^*
\left(I_s\otimes Q(z)-M\right)^{-1}
\ell,
\end{equation}
where
\[
D=-D^*\in\mathbb C^{m\times m},
\qquad
\ell\in\mathbb C^{sr\times m},
\qquad
M=-M^*\in\mathbb C^{sr\times sr},
\]
if and only if there exists a completely positive map
\begin{equation}\label{eq:Ldomain}
\mathcal L:\mathbb C^{rn\times rn}\longrightarrow\mathbb C^{nm\times nm}.
\end{equation}
such that $\mathcal L|_{\Sigma_Q}=\Phi_Q$ and
\begin{equation}\label{eq:KerInclQ}
\begin{split}
\ker\left[
I_m\otimes E_{ij}+\mathcal L(I_r\otimes E_{ij})
\right]_{i,j=1}^n
\subseteq
\ker\biggl[
&B^*(I_m\otimes E_{ij})B\\
&+\mathcal L\left(Q_0^*(I_r\otimes E_{ij})Q_0\right)
\biggr]_{i,j=1}^n.
\end{split}
\end{equation}
\end{theorem}

\begin{proof}
Suppose first that $f(A)=B$ for some
$f\in\mathcal P^{m\times m}_{\rm nc}(\mathbb H_Q)$. Let
$k\in\mathbb N$ and let $T\in\mathbb C^{nk\times nk}$ be invertible. Put
\[
\widetilde H_j=T(\oplus_{\nu=1}^kH_{0,j})T^{-1}.
\]
Since $Q$ is linear,
\begin{equation}\label{eq:Q-similarity}
Q(\widetilde H)
=(I_r\otimes T)(\oplus_{\nu=1}^kQ_0)(I_r\otimes T^{-1}).
\end{equation}
Hence, whenever $\operatorname{Re}Q(\widetilde H)>0$, the tuple
$\widetilde H$ belongs to $\mathbb H_{Q,nk}$. Similarity invariance of nc
rational functions gives
\begin{equation}\label{eq:F-similarity}
f(\widetilde H)
=(I_m\otimes T)(\oplus_{\nu=1}^kB)(I_m\otimes T^{-1}).
\end{equation}
Let $P=T^*T>0$. Multiplying the strict inequality for
$Q(\widetilde H)+Q(\widetilde H)^*$ on the left by $I_r\otimes T^*$ and on the right by $I_r\otimes T$ yields
\[
(I_r\otimes P)(\oplus Q_0)+(
\oplus Q_0^*)(I_r\otimes P)>0.
\]
Likewise,
\[
(I_m\otimes P)(\oplus B)+(
\oplus B^*)(I_m\otimes P)\geq0.
\]
Writing $P=(P_{ij})_{i,j=1}^k$, this says that
\[
\left[(I_r\otimes P_{ij})Q_0+Q_0^*(I_r\otimes P_{ij})\right]_{i,j=1}^k>0
\]
implies
\[
\left[(I_m\otimes P_{ij})B+B^*(I_m\otimes P_{ij})\right]_{i,j=1}^k\geq0.
\]
As usual, a perturbation argument replaces the strict inequality by a
non-strict one. Thus $\operatorname{id}_{\mathbb C^{k\times k}}\otimes\Phi_Q$
is positive for every $k$, and $\Phi_Q$ is completely positive.

Conversely, suppose that a completely positive extension $\mathcal L$
satisfying \eqref{eq:Ldomain} and \eqref{eq:KerInclQ} exists. Choose a
minimal Kraus representation
\[
\mathcal L(Y)=\sum_{\nu=1}^s C_\nu^*YC_\nu,
\qquad
C_\nu\in\mathbb C^{rn\times mn}.
\]
Let $e_1,\ldots,e_n$ be the standard basis of $\mathbb C^n$ and define
\begin{equation}\label{eq:LmatrixQ}
L=\begin{pmatrix}
I_m\otimes e_1^*&\cdots&I_m\otimes e_n^*\\
(I_r\otimes e_1^*)C_1&\cdots&(I_r\otimes e_n^*)C_1\\
\vdots&&\vdots\\
(I_r\otimes e_1^*)C_s&\cdots&(I_r\otimes e_n^*)C_s
\end{pmatrix}
\end{equation}
and
\begin{equation}\label{eq:NmatrixQ}
N=\begin{pmatrix}
(I_m\otimes e_1^*)B&\cdots&(I_m\otimes e_n^*)B\\
-(I_r\otimes e_1^*)Q_0C_1&\cdots&-(I_r\otimes e_n^*)Q_0C_1\\
\vdots&&\vdots\\
-(I_r\otimes e_1^*)Q_0C_s&\cdots&-(I_r\otimes e_n^*)Q_0C_s
\end{pmatrix}.
\end{equation}
The kernel inclusion implies, exactly as in Lemma~\ref{L:KNLequations}, the
existence of a skew-Hermitian matrix
\begin{equation}\label{eq:SQ}
S=\begin{pmatrix}D&\ell^*\\-\ell&-M\end{pmatrix},
\qquad S=-S^*,
\end{equation}
such that $SL=N$. The first block row and the remaining block rows give,
respectively,
\begin{align}
D\otimes I_n+(\ell^*\otimes I_n)
\begin{pmatrix}C_1\\ \vdots\\ C_s\end{pmatrix}&=B,
\label{eq:F0-first-Q}\\
-\ell\otimes I_n-(M\otimes I_n)
\begin{pmatrix}C_1\\ \vdots\\ C_s\end{pmatrix}
&=-(I_s\otimes Q_0)
\begin{pmatrix}C_1\\ \vdots\\ C_s\end{pmatrix}.
\label{eq:C-equation-Q}
\end{align}
Consequently
\[
\begin{pmatrix}C_1\\ \vdots\\ C_s\end{pmatrix}
=
\left(I_s\otimes Q_0-M\otimes I_n\right)^{-1}(\ell\otimes I_n),
\]
and hence $f$ in \eqref{eq:F-realization} satisfies $f(A)=B$.

If $H\in\mathbb H_{Q,N}$ and
\[
R(H)=I_s\otimes Q(H)-M\otimes I_N,
\]
then
\[
R(H)+R(H)^*=I_s\otimes(Q(H)+Q(H)^*)>0.
\]
Thus $R(H)$ is invertible and
\[
f(H)+f(H)^*
=(\ell^*\otimes I_N)R(H)^{-*}(R(H)+R(H)^*)R(H)^{-1}(\ell\otimes I_N)\geq0.
\]
Finally, $Q(-H)=-Q(H)$ and $Q(H)^*=Q(H^*)$, so from
$D^*=-D$ and $M^*=-M$ one obtains
\[
f(-H)^*=-f(H^*).
\]
Thus $f\in\mathcal{PO}^{m\times m}_{\rm nc}(\mathbb H_Q)$.

For the converse, suppose that an interpolant of the form
\eqref{eq:F-realization} exists. Define
\begin{equation}\label{cLform2}
\cL(Y)=  (\ell^* \otimes I_{n})(I_s \otimes Q_0 -M \otimes I_{n})^{-*}(I_s\otimes Y) (I_s \otimes Q_0  -M \otimes I_{n})^{-1}(\ell \otimes I_{n}).
\end{equation}
Clearly, $\mathcal L$ is completely positive. Next, let us check that ${\mathcal L}|_{\Sigma_Q} = \Phi$. Observe first that if \eqref{eq:F-realization}  holds, then \begin{equation}\label{fmulti3} 
B=f(A)= D\otimes I_n +(\ell^* \otimes I_n) (I_s \otimes Q_0- M\otimes I_n)^{-1}(\ell \otimes I_n) . 
 \end{equation} Now
$$
\cL(\left((I_r\otimes X) Q_0 + Q_0^* (I_r\otimes X)\right))$$ $$=(\ell^* \otimes I_n) (I_{s} \otimes Q_0^* +M \otimes I_n)^{-1}(I_s \otimes [(I_r\otimes X) Q_0 +$$ $$ Q_0^* (I_r\otimes X)])
(I_{s} \otimes Q_0-M \otimes I_n)^{-1} (\ell\otimes I_n) .$$
Inserting $0=(I_{rs}\otimes X)(M\otimes I_n) - (M\otimes I_n)(I_{rs}\otimes X)$ in the middle, and using that some factors cancel with their inverses, we obtain that the above equals
$$ (\ell^* \otimes I_n) (I_{rs} \otimes X) (I_{s} \otimes Q_0^* +M \otimes I_n)^{-1} (\ell \otimes I_n) + $$ $$(\ell^* \otimes I_n)  (I_{s} \otimes Q_0 -M \otimes I_n)^{-1} (I_{rs} \otimes X)(\ell \otimes I_n)=$$ 
$$ (I_m\otimes X)f(A) + f(A)^* (I_m\otimes X) = (I_m\otimes X)B + B^* (I_m\otimes X),$$ where we used \eqref{fmulti3}, $I_{rs}\ell=\ell I_m$ and $D=-D^*$. 
This gives
$\mathcal L|_{\Sigma_Q}=\Phi_Q$. 

Moreover, for the Kraus operators arising
from $\cL$, the corresponding matrices $L$ and $N$
satisfy $SL=N$ for the skew-Hermitian matrix in \eqref{eq:SQ}. Hence
$\ker L\subseteq\ker N$, and therefore $\ker L^*L\subseteq\ker N^*N$.
But
\[
L^*L=
\left[I_m\otimes E_{ij}+\mathcal L(I_r\otimes E_{ij})\right]_{i,j=1}^n
\]
and
\[
N^*N=
\left[B^*(I_m\otimes E_{ij})B+
\mathcal L\left(Q_0^*(I_r\otimes E_{ij})Q_0\right)\right]_{i,j=1}^n,
\]
which proves \eqref{eq:KerInclQ}.
\end{proof}

\begin{remark}
    When $Q(z)$ has a direct sum decomposition $Q(z)=\oplus_j Q_j(z)$, then instead of $I_s\otimes Q(z)$ in \eqref{eq:F-realization} one may use a more refined version $\oplus_j (I_{s_j} \otimes Q_j(z))$, as was done,  for example, in Theorem 3.1 in \cite{VW}.
\end{remark}

In the next section we will specify the above result for the poly-half-plane (where $Q(z) = \diag (z_i)_{i=1}^d$). In that case, we can give a function theoretic characterization of which functions $f(z)$ have a realization \eqref{eq:F-realization}; see Theorem \ref{POd}.

\section{The noncommutative poly-half-plane}
\label{sec:polyhalfplane-interpolation}

We now specialize the previous section to the standard
noncommutative right poly-half-plane, where we also obtain some additional results. Take
\[
r=d,
\qquad
L_j=E_{jj}\in\mathbb C^{d\times d},
\qquad j=1,\ldots,d.
\]
Then
\begin{equation}\label{eq:Q-poly}
Q(A)=\sum_{j=1}^dE_{jj}\otimes A_j
=\operatorname{Diag}(A_1,\ldots,A_d).
\end{equation}
Consequently
\[
\mathbb H_{Q,n}=\mathbb{RHP}_n
=\left\{(A_1,\ldots,A_d):\operatorname{Re}A_j>0,\ j=1,\ldots,d\right\}.
\]
Thus the general classes above become the classes
$\mathcal P_{{\rm nc},d}^{m\times m}$ and
$\mathcal{PO}_{{\rm nc},d}^{m\times m}$ used below.

We first record the realization theorem that identifies the relevant
proper part.  This is the result that will allow us to express the final
interpolation theorem in terms of regularity after inversion.

\begin{theorem}\label{POd}
Let $f(z_1,\ldots,z_d)$ be an $m\times m$ nc rational matrix function
in the noncommutative variables $z_1,\ldots,z_d$. The following are
equivalent:
\begin{itemize}
\item[(i)] $f\in\mathcal{PO}^{m\times m}_{{\rm nc},d}$ and
$f(z_1^{-1},\ldots,z_d^{-1})$ is regular at $0$;
\item[(ii)] $f$ has a realization
\begin{equation}\label{fmulti}
f(z)=D+\ell^*
\left(\operatorname{Diag}(z_jI_{r_j})_{j=1}^d-M\right)^{-1}\ell,
\end{equation}
where $r_1,\ldots,r_d\in\mathbb N_0$,
$D=-D^*\in\mathbb C^{m\times m}$,
$\ell\in\mathbb C^{(r_1+\cdots+r_d)\times m}$, and
$M=-M^*\in\mathbb C^{(r_1+\cdots+r_d)\times(r_1+\cdots+r_d)}$.
\end{itemize}
Moreover, if $f$ is real, then $D,\ell,M$ may be chosen real.
\end{theorem}

\begin{proof} $(i)\Rightarrow(ii)$.
Put
\[
g(w_1,\ldots,w_d)
=
f(w_1^{-1},\ldots,w_d^{-1}).
\]
By assumption, $g$ is regular at $0$.
 Define
\[
\Phi(w)=i\,g(w).
\]
We claim that $\Phi$ is matrix-selfadjoint on the noncommutative
imaginary axis.  Indeed, suppose that
\[
W=(W_1,\ldots,W_d),
\qquad
W_j^*=-W_j,
\]
and that the matrices $W_j$ are invertible.  Then
\[
(W_j^{-1})^*=-W_j^{-1},
\]
so that, by the oddness of $f$,
\[
f(W_1^{-1},\ldots,W_d^{-1})
=
-f(W_1^{-1},\ldots,W_d^{-1})^*.
\]
Consequently,
\[
\Phi(W)^*
=
-i\,g(W)^*
=
i\,g(W)
=
\Phi(W).
\]
Since both sides are rational functions and the invertible
skew-Hermitian tuples form a uniqueness set, this identity extends
to all points on the nc imaginary axis at which $\Phi$ is defined.

We can therefore apply Theorem~7.1 of \cite{Alpay2006}.
Thus $\Phi$ has a minimal Givone--Roesser realization
\begin{equation}\label{eq:Phi-GR}
\Phi(w)
=
D_0+
C\bigl(I-\Delta(w)A\bigr)^{-1}\Delta(w)B,
\end{equation}
where
\[
\Delta(w)
=
\operatorname{Diag}(w_iI_{r_i})_{i=1}^d,
\qquad
D_0=D_0^*,
\]
and there is an invertible structured Hermitian matrix
\[
H=\operatorname{Diag}(H_1,\ldots,H_d)
\]
such that
\begin{equation}\label{eq:selfadjoint-relations}
A^*H+HA=0,
\qquad
C=iB^*H.
\end{equation}

For completeness, the connection with the $J$-unitary theory in
\cite{Alpay2006} is as follows.  Theorem~7.1 is obtained by associating
to $\Phi$ the block function
\[
F(w)
=
\begin{pmatrix}
I_m&i\Phi(w)\\
0&I_m
\end{pmatrix},
\qquad
J_1=
\begin{pmatrix}
0&I_m\\
I_m&0
\end{pmatrix}.
\]
The matrix-selfadjoint property of $\Phi$ on the nc imaginary axis
is equivalent to the $J_1$-unitarity of $F$ there.

We now use the positive-real hypothesis.  If
\[
W=(W_1,\ldots,W_d)\in\mathbb{RHP}_n,
\]
then each $W_j$ is invertible and
\[
\operatorname{Re}W_j^{-1}
=
W_j^{-*}(\operatorname{Re}W_j)W_j^{-1}>0.
\]
Hence
\[
(W_1^{-1},\ldots,W_d^{-1})\in\mathbb{RHP}_n.
\]
Since $f$ is positive real,
\[
\operatorname{Re}g(W)\geq0.
\]
As $\Phi=ig$, this says
\[
\operatorname{Im}\Phi(W)
=
\frac{\Phi(W)-\Phi(W)^*}{2i}
=
\operatorname{Re}g(W)
\geq0.
\]

A direct calculation gives
\[
J_1-F(W)J_1F(W)^*
=
\begin{pmatrix}
2\operatorname{Im}\Phi(W)&0\\
0&0
\end{pmatrix}
\geq0.
\]
Thus the $J_1$-unitary function $F$ occurring in the proof of
Theorem~7.1 is in fact $J_1$-inner on the nc right poly-half-plane.
By Theorem~6.1 of \cite{Alpay2006}, its associated structured
Hermitian matrix is therefore strictly positive:
\[
H>0.
\]

Because $H$ is block diagonal, so is $H^{1/2}$.  Performing the
block-diagonal similarity
\[
A\longmapsto H^{1/2}AH^{-1/2},
\qquad
B\longmapsto H^{1/2}B,
\qquad
C\longmapsto CH^{-1/2},
\]
we may therefore assume without loss of generality that
\[
H=I_r.
\]
The relations \eqref{eq:selfadjoint-relations} then become
\[
A^*=-A,
\qquad
C=iB^*.
\]
Consequently \eqref{eq:Phi-GR} takes the form
\[
\Phi(w)
=
D_0+iB^*
\bigl(I-\Delta(w)A\bigr)^{-1}
\Delta(w)B,
\qquad D_0=D_0^*.
\]
Since $g=-i\Phi$, we obtain
\[
g(w)
=
D+B^*
\bigl(I-\Delta(w)A\bigr)^{-1}
\Delta(w)B,
\qquad
D=-iD_0,
\]
and hence
\[
D^*=-D.
\]

Now substitute $w_i=z_i^{-1}$.  Since
\[
\Delta(z^{-1})=\Delta(z)^{-1},
\]
we have
\begin{align*}
\bigl(I-\Delta(z^{-1})A\bigr)^{-1}\Delta(z^{-1})
&=
\bigl(I-\Delta(z)^{-1}A\bigr)^{-1}\Delta(z)^{-1}\\
&=
\bigl(\Delta(z)-A\bigr)^{-1}.
\end{align*}
It follows that
\[
f(z)
=
D+B^*\bigl(\Delta(z)-A\bigr)^{-1}B.
\]
Setting
\[
M=A,
\qquad
\ell=B,
\]
we have
\[
M^*=-M
\]
and therefore
\[
f(z)
=
D+\ell^*
\left(
\operatorname{Diag}(z_iI_{r_i})_{i=1}^d-M
\right)^{-1}\ell,
\]
which is \eqref{fmulti}.

We next prove $(ii)\Rightarrow(i)$.  Write
\[
\Delta(Z)=\operatorname{Diag}(I_{r_i}\otimes Z_i)_{i=1}^d
\]
at matrix level.  If
\[
Z=(Z_1,\ldots,Z_d)\in\mathbb{RHP}_n,
\]
then
\[
\operatorname{Re}\Delta(Z)>0.
\]
Since $M^*=-M$,
\[
\operatorname{Re}\bigl(\Delta(Z)-M\otimes I_n\bigr)
=
\operatorname{Re}\Delta(Z)>0.
\]
In particular, $\Delta(Z)-M\otimes I_n$ is invertible.  Thus
\eqref{fmulti} is defined on every nc right poly-half-plane.

For every invertible matrix $X$,
\[
\operatorname{Re}X^{-1}
=
X^{-*}(\operatorname{Re}X)X^{-1}.
\]
Applying this to
\[
X=\Delta(Z)-M\otimes I_n
\]
gives
\[
\operatorname{Re}X^{-1}\geq0.
\]
Since $D^*=-D$, it follows that
\[
\operatorname{Re}f(Z)
=
(\ell^*\otimes I_n)
\operatorname{Re}(X^{-1})
(\ell\otimes I_n)
\geq0.
\]
Hence $f\in\cP_{nc,d}^{m\times m}$.

The oddness is also immediate.  Indeed, putting
\[
X=\Delta(Z)-M\otimes I_n,
\]
we have
\[
\Delta(-Z^*)-M\otimes I_n
=
-X^*.
\]
Therefore
\begin{align*}
-f(-Z^*)^*
&=
-D^*
-
\left[
(\ell^*\otimes I_n)(-X^*)^{-1}(\ell\otimes I_n)
\right]^*\\
&=
D+
(\ell^*\otimes I_n)X^{-1}(\ell\otimes I_n)\\
&=
f(Z).
\end{align*}
Thus
\[
f\in\cO_{nc,d}^{m\times m},
\]
and consequently
\[
f\in\PO_{nc,d}^{m\times m}.
\]

It remains to verify regularity at $0$ after inversion.  From
\eqref{fmulti},
\[
f(w_1^{-1},\ldots,w_d^{-1})
=
D+
\ell^*
\bigl(\Delta(w)^{-1}-M\bigr)^{-1}\ell.
\]
Using
\[
\Delta(w)^{-1}-M
=
\Delta(w)^{-1}\bigl(I-\Delta(w)M\bigr),
\]
we obtain
\[
\bigl(\Delta(w)^{-1}-M\bigr)^{-1}
=
\bigl(I-\Delta(w)M\bigr)^{-1}\Delta(w).
\]
Hence
\[
f(w_1^{-1},\ldots,w_d^{-1})
=
D+
\ell^*
\bigl(I-\Delta(w)M\bigr)^{-1}
\Delta(w)\ell.
\]
The expression on the right is defined at $w=0$, and therefore
$f(z_1^{-1},\ldots,z_d^{-1})$ is regular at $0$.

Finally, suppose that $f$ is real.  The preceding construction can
be carried out over $\mathbb R$: a minimal real GR-realization may
be chosen for the real rational function, and the realization
identities and the normalization by the positive structured
Hermitian matrix preserve the real structure.  Thus $D,\ell$, and
$M$ may be chosen real.  Since they also satisfy
\[
D^*=-D,\qquad M^*=-M,
\]
in the real case they are skew-symmetric:
\[
D^T=-D,\qquad M^T=-M.
\]
Conversely, if $D,\ell$, and $M$ in \eqref{fmulti} are real, then
the rational expression \eqref{fmulti} has real coefficients, and
hence $f\in\PRO_{nc,d}^{m\times m}$.
\end{proof}

We can now state the matrix-point interpolation theorem in a form that
retains the regularity condition from Theorem~\ref{POd}.

\begin{theorem}\label{multi}
Let
\[
A=(A_1,\ldots,A_d)\in\mathbb{RHP}_n,
\qquad
B\in\mathbb C^{nm\times nm},
\]
and put
\[
\mathcal A=\operatorname{Diag}(A_1,\ldots,A_d).
\]
Define
\begin{equation}\label{eq:Sigma-poly}
\Sigma=
\left\{
(I_d\otimes X)\mathcal A+\mathcal A^*(I_d\otimes X):
X\in\mathbb C^{n\times n}
\right\}
\end{equation}
and
\begin{equation}\label{eq:Phi-poly}
\Phi\left((I_d\otimes X)\mathcal A+\mathcal A^*(I_d\otimes X)\right)
=(I_m\otimes X)B+B^*(I_m\otimes X).
\end{equation}

If there exists $f\in\mathcal P^{m\times m}_{{\rm nc},d}$ with $f(A)=B$,
then $\Phi$ is completely positive.

Moreover, there exists an interpolant
\[
f\in\mathcal{PO}^{m\times m}_{{\rm nc},d}
\]
such that
\begin{equation}\label{eq:regular-interpolant}
f(A)=B,
\qquad
f(z_1^{-1},\ldots,z_d^{-1})\text{ is regular at }0,
\end{equation}
if and only if there exists a completely positive map
\begin{equation}\label{eq:Ldomain2}
\mathcal L:\mathbb C^{nd\times nd}\longrightarrow\mathbb C^{nm\times nm}.
\end{equation}
such that $\mathcal L|_\Sigma=\Phi$ and
\begin{equation}\label{KerIncl2}
\begin{split}
\ker
\left[
I_m\otimes E_{ij}+\mathcal L(I_d\otimes E_{ij})
\right]_{i,j=1}^n
\subseteq
\ker\biggl[
&B^*(I_m\otimes E_{ij})B\\
&+\mathcal L\left(\operatorname{Diag}(A_s^*E_{ij}A_s)_{s=1}^d\right)
\biggr]_{i,j=1}^n.
\end{split}
\end{equation}
\end{theorem}

\begin{proof}
The first assertion follows immediately from the first part of
Theorem~\ref{thm:Q-interpolation} with the choice
$r=d$ and $L_j=E_{jj}$.

For the second assertion, with the same choice of pencil we have
\[
Q(A)=\mathcal A,
\qquad
\Sigma_Q=\Sigma,
\qquad
\Phi_Q=\Phi,
\]
and
\[
Q(A)^*(I_d\otimes E_{ij})Q(A)
=\operatorname{Diag}(A_s^*E_{ij}A_s)_{s=1}^d.
\]
Thus the kernel condition \eqref{eq:KerInclQ} is exactly
\eqref{KerIncl2}. By Theorem~\ref{thm:Q-interpolation}, the extension
condition is equivalent to the existence of an interpolant of the proper
realization type
\[
f(z)=D+\ell^*
\left(\operatorname{Diag}(z_jI_{r_j})_{j=1}^d-M\right)^{-1}\ell.
\]
By Theorem~\ref{POd}, this is equivalent to the existence of an
interpolant satisfying \eqref{eq:regular-interpolant}.
\end{proof}

Theorem \ref{multi} suggests the following algorithm for finding an interpolant.
We formulate the algorithm for the poly-half-plane, but a similar algorithm can be formulated for the more general domains.

\vfill

\begin{algorithm}
\caption{Finding a nc rational matrix interpolant}
\label{alg:matrix-interpolant}
\begin{algorithmic}[1]
\Procedure{FindInterpolant}{$(A_1,\ldots,A_d),B$}
\State Set
\[
\mathcal A=\operatorname{Diag}(A_1,\ldots,A_d)
\]
and define $\Sigma$ by \eqref{eq:Sigma-poly}.
\State Define $\Phi$ by \eqref{eq:Phi-poly}.
\If{$\Phi:\Sigma\to\mathbb C^{nm\times nm}$ is not completely positive}
    \State \Return ``No interpolant in $\mathcal P^{m\times m}_{{\rm nc},d}$ exists.''
\Else
    \State Search for a completely positive Arveson extension
    \[
    \mathcal L:\mathbb C^{nd\times nd}\to\mathbb C^{nm\times nm},
    \qquad \mathcal L|_\Sigma=\Phi,
    \]
    satisfying \eqref{KerIncl2}.
    \If{such an extension $\mathcal L$ exists}
        \State Find a minimal Kraus representation
        \[
        \mathcal L(Y)=\sum_{\nu=1}^s C_\nu^*YC_\nu,
        \qquad C_\nu\in\mathbb C^{nd\times nm}.
        \]
        \State Form $L$ and $N$ as in \eqref{eq:LmatrixQ} and
        \eqref{eq:NmatrixQ}, with $r=d$ and $Q_0=\mathcal A$.
        \State Solve
        \[
        SL=N,\qquad S=-S^*,
        \]
        and write
        \[
        S=\begin{pmatrix}D&\ell^*\\-\ell&-M\end{pmatrix}.
        \]
        \State Define
        \[
        f(z)=D+\ell^*
        \left(I_s\otimes\operatorname{Diag}(z_1,\ldots,z_d)-M\right)^{-1}\ell.
        \]
        \State \Return $f$.
        \Comment{$f(A)=B$ and $f(z^{-1})$ is regular at $0$.}
    \Else
        \State Search for matrices
        \[
        K_1,\ldots,K_d\in\mathbb C^{m\times m},
        \qquad K_j\succeq0,
        \]
     maximizing
     $\sum_{j=1}^d\operatorname{tr}K_j,
        $
        subject to complete positivity of the reduced map corresponding to
        \[
        \widehat B=B-\sum_{j=1}^dK_j\otimes A_j.
        \]
        \State Apply \textsc{FindInterpolant} to
        $(A_1,\ldots,A_d),\widehat B$.
        \If{this produces an interpolant $\widehat f$}
            \State Define
            \[
            f(z)=\sum_{j=1}^d z_jK_j+\widehat f(z).
            \]
            \State \Return $f$. \Comment{$f(A)=B$ and $f \in \mathcal{PO}^{m\times m}_{{\rm nc},d}$.}
        \Else
            \State \Return ``Inconclusive.''
            \Comment{A possibility when $d>1$.}
        \EndIf
    \EndIf
\EndIf
\EndProcedure
\end{algorithmic}
\end{algorithm}

\eject

Some remarks are in order.
\begin{itemize}
    \item For Step 4: Checking whether $\Phi$ is completely positive is equivalent (as the operator system contains a positive definite element) to checking whether it has an Arveson extension. The latter is a SemiDefinite Program (SDP) and can thus be checked numerically. Typically one would skip Step 3 and see what Step 7 brings.
    \item For Step 7: When $P$ and $Q$ are positive semidefinite, the inclusion ${\rm ker}P \subseteq {\rm ker Q}$ is equivalent to the existence of an $x>0$ such that $P-xQ \ge 0$. Thus, for a fixed $x>0$ Step 7 can be implemented as an SDP. For the subsequent steps it is probably advantageous to maximize $x$, as the Kraus representation of the map ${\mathcal L} $ will be simpler (in fact, it seems that $\| S\| = \frac{1}{\sqrt{x_{\max}}}$). Once $x_{\max}$ is found one may minimize the trace of the Choi matrix of ${\mathcal L}$ to (heuristically) minimize the rank of the Choi matrix.
    \item As the operator system has a diagonal form we will typically find a block diagonal Choi matrix for ${\mathcal L}$. We may now change the off-diagonal matrix of the Choi matrix such that the rank of the Choi matrix will equal the maximal rank of the diagonal blocks (again to find simpler Kraus representation).
    \item When the above algorithm fails to produce a solution, one may try to find one by using that the class $\PO^{m\times m}_{nc,d}$ is invariant under replacing $z_i$ by $\frac{1}{z_i}$, $i=1,\ldots, d$, and also if $f \in \PO^{m\times m}_{nc,d}$ is invertible, then $f^{-1} \in\PO^{m\times m}_{nc,d}$. Thus, one can rerun the algorithm with some $A_i$s or $B$ replaced by their inverse, and then when a solution is found reverse the corresponding changes. 
\end{itemize}

Let us revisit Example \ref{ex2} from the above point of view.
\begin{example}\label{Ex:Lp}
Let $A={\rm Diag} (1,2,3)$ and $B = {\rm Diag} (2,\frac52,\frac{10}{3})$. 
The completely positive map ${\mathcal L}$ is $X \mapsto H\circ X$, where $\circ$ denotes the Schur product and $H$ is the Pick matrix
$$ P = \begin{pmatrix} 2 & \frac32 & \frac43 \cr \frac32 & \frac54 & \frac76 \cr \frac43 & \frac76 & \frac{10}{9} \end{pmatrix}. $$
Checking condition \eqref{KerIncl} we get 
$$  [E_{ij}+{\mathcal L}(E_{ij})]_{i,j=1}^3=\begin{pmatrix} 3 & 0 & 0 & 0 & \frac52 & 0 & 0 & 0 & \frac73 \cr 
0 & 0 & 0 & 0 & 0 & 0 & 0 & 0 & 0 \cr 
0 & 0 & 0 & 0 &0 & 0 & 0 & 0 & 0 \cr 
0 & 0 & 0 & 0 &0 & 0 & 0 & 0 & 0 \cr 
\frac52 & 0 & 0 & 0 & \frac94 & 0 & 0 & 0 & \frac{13}{6} \cr
0 & 0 & 0 & 0 &0 & 0 & 0 & 0 & 0 \cr 
0 & 0 & 0 & 0 &0 & 0 & 0 & 0 & 0 \cr 
0 & 0 & 0 & 0 &0 & 0 & 0 & 0 & 0 \cr  
\frac73 & 0 & 0 & 0 &\frac{13}{6} & 0 & 0 & 0 & \frac{19}{9} \end{pmatrix}, $$ and
$$   [ B^*E_{ij}B + {\mathcal L} (A^* E_{ij} A)_{i,j}]_{i,j=1}^3  = \begin{pmatrix} 6 & 0 & 0 & 0 & 8 & 0 & 0 & 0 & 10\frac23 \cr 
0 & 0 & 0 & 0 & 0 & 0 & 0 & 0 & 0 \cr 
0 & 0 & 0 & 0 &0 & 0 & 0 & 0 & 0 \cr 
0 & 0 & 0 & 0 &0 & 0 & 0 & 0 & 0 \cr 
8 & 0 & 0 & 0 & 11\frac14 & 0 & 0 & 0 & 15\frac13 \cr
0 & 0 & 0 & 0 &0 & 0 & 0 & 0 & 0 \cr 
0 & 0 & 0 & 0 &0 & 0 & 0 & 0 & 0 \cr 
0 & 0 & 0 & 0 &0 & 0 & 0 & 0 & 0 \cr  
10\frac23 & 0 & 0 & 0 &15\frac13 & 0 & 0 & 0 & 21\frac{1}{9} \end{pmatrix}. $$
We have that 
$$ \begin{pmatrix} 1 & 0 & 0  &0  & -4 & 0 & 0 & 3 \end{pmatrix}^T$$ belongs to the kernel of the first matrix, but not to the kernel of the second one. Thus \eqref{KerIncl} is violated.

Observe that $$  [{\mathcal L}(E_{i,j})]_{i,j=1}^3=\begin{pmatrix} 2 & 0 & 0 & 0 & \frac32 & 0 & 0 & 0 & \frac43 \cr 
0 & 0 & 0 & 0 & 0 & 0 & 0 & 0 & 0 \cr 
0 & 0 & 0 & 0 &0 & 0 & 0 & 0 & 0 \cr 
0 & 0 & 0 & 0 &0 & 0 & 0 & 0 & 0 \cr 
\frac32 & 0 & 0 & 0 & \frac54 & 0 & 0 & 0 & \frac{7}{6} \cr
0 & 0 & 0 & 0 &0 & 0 & 0 & 0 & 0 \cr 
0 & 0 & 0 & 0 &0 & 0 & 0 & 0 & 0 \cr 
0 & 0 & 0 & 0 &0 & 0 & 0 & 0 & 0 \cr  
\frac43 & 0 & 0 & 0 &\frac{7}{6} & 0 & 0 & 0 & \frac{10}{9} \end{pmatrix} $$ has $w$ in its range.
Subtracting the maximal multiple of $ww^*$ we find $a_0=1$, and obtain with $\hat{\mathcal L}= {\mathcal L} - {\rm id}$ the equation
$$ S \begin{pmatrix} 1 & 0 & 0 & 0 & 1 & 0 & 0 & 0 & 1 \cr 1 & 0 & 0 & 0 & \frac12 & 0 & 0 & 0 & \frac13 \end{pmatrix}
= \begin{pmatrix}  1 & 0 & 0 & 0 & \frac12 & 0 & 0 & 0 & \frac13 \cr -1 & 0 & 0 & 0 & -1 & 0 & 0 & 0 & -1 \end{pmatrix}.
$$ Solving yields 
$$
S =
\begin{pmatrix}
   0 & 1 \cr -1 & 0 
\end{pmatrix},$$ giving the solution $\frac{1}{z}$. Thus the final solution is $f(z)=z+\frac1z$.
\end{example} 

Next, we present a two-variable example.

\begin{example}\rm

Let $$ A_{1} = \begin{pmatrix}
    1 & 0 \cr 0 & 2
\end{pmatrix}, A_{2} = \begin{pmatrix}
    2 & 1 \cr 1 & 2
\end{pmatrix}, B=\begin{pmatrix} 1 & -2/13\cr  0 & 11/13 \end{pmatrix}.$$
For $X=\begin{pmatrix}
    p & q \cr r & s
\end{pmatrix}$, we get that $\mathcal L$ should act
\begin{equation}\label{condL} \begin{pmatrix} 2p & 3q &           0 &           0\cr 
3r & 4s &           0&           0 \cr
 0 &   0 & 4p + q + r & p + 4q + s\cr 
  0 &   0&  p + 4r + s & q + r + 4s \end{pmatrix} \longmapsto \begin{pmatrix}                 2p &            24q/13 - 2p/13 \cr
24r/13 - 2p/13 & 22s/13 - 2r/13 - 2q/13 \end{pmatrix}.\end{equation}
To find an Arveson extension ${\mathcal L}: {\mathbb C}^{4\times 4} \to {\mathbb C}^{2\times 2},$ we can find its Choi matrix
$$ L=[L_{ij}]_{i,j=1}^4:= [{\mathcal L}(E_{ij})]_{i,j=1}^4 \in {\mathbb C}^{8\times 8}\cap {\rm PSD}.$$
Here $E_{ij} \in {\mathbb C}^{4\times 4}$ has a 1 in entry $(i,j)$ and zeros elsewhere, and $L_{ij} \in {\mathbb C}^{2\times 2}.$
Condition \eqref{condL} implies the following linear constraints on $L$:
\begin{equation}\label{11}
    2L_{11}+4L_{33}
+L_{34}+L_{43}=\begin{pmatrix}                 2 &            -2/13 \cr
- 2/13 & 0 \end{pmatrix},\end{equation}
\begin{equation}\label{12}
    3L_{12}+L_{33}
+4L_{34}+L_{44}=\begin{pmatrix}                 0 &            24/13 \cr
0 & -2/13 \end{pmatrix},\end{equation}
\begin{equation}\label{13}
    3L_{21}+L_{33}
+4L_{43}+L_{44}=\begin{pmatrix}                 0 &  0\cr           24/13 & -2/13 \end{pmatrix},\end{equation}
\begin{equation}\label{14}
    4L_{22}+L_{34}
+L_{43}+4L_{44}=\begin{pmatrix}                 0 &  0\cr           0 & 22/13 \end{pmatrix},\end{equation}
Note that there are no conditions on $L_{ij}, L_{ji}$ for $(i,j)\in \{ (1,3),(1,4),(2,3),(2,4)\}$, and the existence of an $L\in {\rm PSD}$ implies that $[L_{ij}]_{i,j=1}^2 \in {\rm PSD}$ and $[L_{ij}]_{i,j=3}^4 \in {\rm PSD}$, so that one may choose $L_{ij}=0= L_{ji}$ for $(i,j)\in \{ (1,3),(1,4),(2,3),(2,4)\}$.

When we perform in Step 7 with $x_{\max}=\frac13$ an SDP minimizing trace among all PSD $L$ satisfying \eqref{11}-\eqref{14}, we find a numerical solution that leads to
$$ L= \begin{pmatrix}  vv^T & 0 \cr 0 &  yy^T \end{pmatrix}, $$
where 
$$ v^T=\begin{pmatrix}  1 & -\frac{1}{13} & 0 & \frac{8}{13} \end{pmatrix} , \ y^T= \begin{pmatrix} 0 & \frac{1}{13} & 0 & -\frac{3}{13}\end{pmatrix}.$$

To get a rank 1 solution, which corresponds to a Kraus representation with 1 term, we can define 
$$\tilde{L}= \begin{pmatrix}  vv^T &  vy^T \cr  yv^T &  yy^T \end{pmatrix}, $$
which is PSD and satisfies \eqref{11}-\eqref{14}. 
We obtain that 
$$ C_1=\begin{pmatrix} 1 & -\frac{1}{13} \cr 0 & \frac{8}{13} \cr 0 & \frac{1}{13} \cr 0 & -\frac{3}{13}, \end{pmatrix}$$
and we obtain the equation
$$ S \begin{pmatrix} 1 & 0 & 0 & 1 \cr 1 & -\frac{1}{13} & 0 & \frac{8}{13}
\cr 0 & \frac{1}{13} & 0 & -\frac{3}{13}, \end{pmatrix}= 
\begin{pmatrix} 1 & -\frac{2}{13} & 0 & \frac{11}{13} \cr -1 & \frac{1}{13} & 0 & -\frac{16}{13}
\cr 0 & \frac{1}{13} & 0 & \frac{5}{13}, \end{pmatrix} $$

Computing the solution we find 
$$ S=\begin{pmatrix}
    0 & 1 & -1 \cr -1 &  0 & 1 \cr 1 & -1 & 0
\end{pmatrix}.$$
Thus
\begin{equation}\label{sol}
D=0, \ell =  \begin{pmatrix}
    1 \cr -1 
\end{pmatrix} , M = \begin{pmatrix}
    0 & -1 \cr 1 & 0
\end{pmatrix}, F(z_1,z_2) = \ell^* \left[\begin{pmatrix}
    z_1 & 0 \cr 0 & z_2
\end{pmatrix}-M \right]^{-1} \ell. \end{equation}

\end{example}

\subsection*{Acknowledgements}
This work is based on research supported in part by the National Research Foundation of South Africa (NRF, Grant Number 127364) and the DSI-NRF Centre of Excellence in Mathematical and Statistical Sciences (CoE-MaSS). Any opinion, finding and conclusion or recommendation expressed in this material is that of the authors and the NRF and CoE-MaSS do not accept any liability in this regard. 
VV thanks the Milken Families Foundation Chair in Mathematics
for partial support of his research. In addition, HW is partially supported by United States National Science Foundation (NSF) grant DMS-2348720.
Vinnikov and Woerdeman wish to
thank the Banff International Research Station for the opportunity to work on this project
during the research in teams program 24rit023.

\subsection*{AI disclosure} (i) GPT-5.6 Sol was used to write a first draft of Subsection 2.1, which recalls various definitions and standard results. Corrections had to be made to this first draft; (ii) Based on an outline provided by the authors (including specific theorems and references), GPT-5.6 Sol was used to fill in the remaining details of the proof of Theorem 5.1; (iii) When the decision was made to reverse the order of Sections 4 and 5, GPT-5.6 Sol was helpful in drafting a first round of this reversal; (iv) GPT-5.6 Sol was used to perform some more trivial tasks such as converting references to BibTeX format and reference search. 



\bibliographystyle{plain}





%
%
\end{document}